\documentclass[11pt]{article}
\usepackage{amsthm, amsmath, amssymb, amsfonts, fancyhdr, verbatim, url}
\usepackage[margin = 1in]{geometry}
\usepackage{graphicx}
\usepackage{setspace}
\usepackage{tikz}
\usepackage[notref, notcite, final]{showkeys}
\usepackage{booktabs}
\usepackage{hyperref}

\usepackage{mathpazo}
\newtheorem*{claim*}{Claim}
\newtheorem*{fact*}{Fact}
\newtheorem{theorem}{Theorem}[section]
\newtheorem*{theorem*}{Theorem}
\newtheorem{proposition}[theorem]{Proposition}
\newtheorem{lemma}[theorem]{Lemma}
\newtheorem*{lemma*}{Lemma}
\newtheorem{corollary}[theorem]{Corollary}
\newtheorem{conjecture}[theorem]{Conjecture}

\theoremstyle{definition}
\newtheorem{definition}[theorem]{Definition}
\newtheorem{example}[theorem]{Example}

\newtheorem*{question*}{Question}

\theoremstyle{remark}

\def\c{\mathcal}

\def\ZZ{\mathbb{Z}}

\def\PP{\mathbb{P}}
\def\EE{\mathbb{E}}

\DeclareMathOperator{\Var}{Var}

\newcommand{\abs}[1]{\left\lvert#1\right\rvert}

\newcommand{\floor}[1]{\left\lfloor #1 \right\rfloor}

\newcommand{\set}[1]{\left\{#1\right\}}
\def\({\left(}
\def\){\right)}

\def\ol#1{\overline{#1}}

\def\ol#1{\overline{#1}}

\def\x{\times}

\def\jbar{\overline \jmath}
\def\kbar{\overline{k}}

\def\<={\Leftarrow}
\def\=>{\Rightarrow}

\title{Sets Characterized by Missing Sums and Differences}
\author{Yufei Zhao \\[6pt]
Massachusetts Institute of Technology\\
{\tt yufeiz@mit.edu}
}

\begin{document}

\maketitle

\begin{abstract}
A more sums than differences (MSTD) set is a finite subset $S$ of the integers such $\abs{S+S} > \abs{S-S}$. We show that the probability that a uniform random subset of $\{0, 1, \dots, n\}$ is an MSTD set approaches some limit $\rho > 4.28 \x 10^{-4}$. This improves the previous result of Martin and O'Bryant that there is a lower limit of at least $2 \x 10^{-7}$. Monte Carlo experiments suggest that $\rho \approx 4.5 \x 10^{-4}$. We present a deterministic algorithm that can compute $\rho$ up to arbitrary precision.

We also describe the structure of a random MSTD set $S \subseteq \set{0,1,\dots,n}$. We formalize the intuition that fringe elements are most significant, while middle elements are nearly unrestricted. For instance, the probability that any ``middle'' element is in $S$ approaches $1/2$ as $n \to \infty$, confirming a conjecture of Miller, Orosz, and Scheinerman.

In general, our results work for any specification on the number of missing sums and the number of missing differences of $S$, with MSTD sets being a special case.
\end{abstract}

%%%%%%%%%%%%%%%%%%%%%%%%%%%%%%%%%%%%%%%%%%%%%%%%%%%%%%%%%%%%%%%%%%%%%%%%%%%%%%%%%%%%%%%%%%%%%%%%

\section{Introduction} \label{sec:intro}

A more sums than differences (MSTD) set is a finite set $S$ of integers with $\abs{S + S} > \abs{S - S}$, where the sum set $S + S$ and the difference set $S - S$ are defined as
\begin{align*}
	S + S & = \{s_1 + s_2 : s_1, s_2 \in S\}, \\
	S - S &= \{s_1 - s_2 : s_1, s_2 \in S\}.
\end{align*}
Since addition is commutative while subtraction is not, two distinct integers $s_1$ and $s_2$ generate one sum but two differences. This suggests that $S+S$ should ``usually'' be smaller than $S - S$. Thus we expect MSTD sets to be rare.

The first example of an MSTD was found by Conway in the 1960's: $\{0,2,3,4,7,11,12,14\}$. The name MSTD was later given by Nathanson
\cite{Na}. MSTD sets have recently become a popular research topic \cite{Hegarty, HM, MO, MOS, Na2, Na, Zhao:bidirectional, Zhao:finite}. For older papers see \cite{HRY, Marica, PF, Roesler, Ru1, Ru2, Ru3}. We refer the reader to \cite{Na2, Na} for the history of the problem.

In this paper, we address the following two questions regarding MSTD sets and their generalizations.
\begin{enumerate}
	\item What is the probability that a random subsets of $\{0, 1, \dots, n\}$ is an MSTD set?
	\item What is the structure of a typical random MSTD subset of $\{0, 1, \dots, n\}$?
\end{enumerate}
The first question was raised by Martin and O'Bryant \cite{MO}. Let $\rho_n$ be the probability that a uniformly chosen random subset of $\{0, 1, \dots, n\}$ is an MSTD set. In \cite{MO} it was shown that $\rho_n \geq 2 \times 10^{-7}$ for all $n \geq 14$. This is a surprising result since it is contrary to our original intuition that MSTD sets should be rare. It is true that $\rho_n = 0$ for $n \leq 13$, and $\rho_n$ is then monotonically increasing at least for $n \leq 26$. From this data, Martin and O'Bryant conjectured that $\rho_n$ approaches some limit and then they estimated this limit using Monte Carlo experiments.

\begin{conjecture}[Martin and O'Bryant \cite{MO}]
As $n \to \infty$, the proportion $\rho_n$ of MSTD sets converges to a limit about $4.5 \x 10^{-4}$.
\end{conjecture}

Previously it was not known whether $\rho_n$ converges. In this paper, we show that $\rho_n$ indeed approaches some limit $\rho$. We also give a deterministic algorithm which can, in principle, compute arbitrarily good lower and upper bounds for $\rho$. 

\begin{theorem} \label{thm:main}
As $n \to \infty$, the proportion $\rho_n$ of MSTD sets converges to a limit $\rho > 4.28 \x 10^{-4}$.
\end{theorem}

Our numerical result is a significant improvement over Martin and O'Bryant's $2 \x 10^{-7}$. Unfortunately, limits of computation prevent us from giving a good upper bound. However, if we were to have unlimited computing power, then our method could give provable bounds for $\rho$ up to any desired precision.

Our proof, like that of Martin and O'Bryant, is non-constructive. As for constructive results, the densest families of MSTDs subsets of $\{0, 1, 2, \dots, n\}$ constructed so far are due to Miller, Orosz, and Scheinerman \cite{MOS} (with density $\Omega(1/n^4)$) and the author \cite{Zhao:bidirectional} (with density $\Theta(1/n)$). No explicit explicit construction with $\Omega(1)$ density is known.

Our method for proving Theorem \ref{thm:main} can easily be adapted to answer other similar questions such as:
\begin{enumerate}
	\item What is the probability that a uniformly random subset $S \subseteq \{0, 1, \dots, n\}$ has more differences than sums, i.e., $\abs{S+S} < \abs{S-S}$?
	\item What is the probability that a uniformly random subset $S \subseteq \{0, 1, \dots, n\}$ has equal number of differences and sums, i.e., $\abs{S+S} = \abs{S-S}$?
	\item What is the probability that a uniformly random subset $S \subseteq \{0, 1, \dots, n\}$ is missing exactly $s$ sums and $d$ differences, i.e., $\abs{S+S} = 2n+1 - s$, $\abs{S-S} = 2n + 1 - d$, where $s$ and $d$ are fixed?
	\item What is the probability that a uniformly random subset $S \subseteq \{0, 1, \dots, n\}$ has exactly $x$ more sums than differences, i.e. $\abs{S+S} - \abs{S-S} = x$, where $x$ is fixed?
\end{enumerate}
As we will show, in each case, as $n \to \infty$, each sequence of probabilities approaches some limit. Furthermore, we have a deterministic algorithm that can give arbitrarily good provable bounds for the limit.

Our general result works for any characterization on the number of missing sums and and the number of missing differences of $S \subseteq \{0, 1, \dots, n\}$, by which we mean the pair
\[
	\lambda(S) = \lambda_n(S) = \( 2n+1 - \abs{S+S}, 2n+1 - \abs{S-S} \).
\]
Let $\Lambda$ denote some (possibly infinite) subset of $\ZZ_{\geq 0} \x \ZZ_{\geq 0}$. We would like to study the collection of subsets $S \subseteq \{0, 1, \dots, n\}$ such that $\lambda(S) \in \Lambda$. For instance, $\Lambda = \{(s,d) : s < d\}$ corresponds to MSTD sets; the one-element set $\Lambda = \{(s,d)\}$ corresponds to question 3 above; $\Lambda = \{(s, d) : d - s = x\}$ corresponds to Question 4 above.

Let
\[
	\rho^\Lambda_n = 2^{-n-1} \abs{\set{ S \subseteq \{0, 1, \dots, n\} : \lambda(S) \in \Lambda }}.
\]
This is the probability that a uniformly random subset of $\{0, 1, \dots, n\}$ characterized by $\Lambda$. We prove the following generalization of Theorem \ref{thm:main}. When $\Lambda$ is the one-element set $\{(s,d)\}$, we abuse notation by writing $\rho^{s,d}$ to mean $\rho^\Lambda$.

\begin{theorem} \label{thm:gen}
For any $\Lambda \subseteq \ZZ_{\geq 0} \x \ZZ_{\geq 0}$, the limit
\[
	\rho^\Lambda = \lim_{n \to \infty} \rho_n^\Lambda
\]
exists. It is positive as long as $\Lambda$ contains as least one element $(s,d)$ where $d$ is even. Furthermore,
\[
	\rho^\Lambda = \sum_{(s,d) \in \Lambda} \rho^{s,d}.
\]
\end{theorem}

Theorem \ref{thm:gen} resolves Conjectures 2 and 19 of Martin and O'Bryant \cite{MO}. Specifically, they conjectured that the probabilities in questions 1--3 above all have limits as $n \to \infty$, and also that $\sum _{s, d} \rho^{s,d} = 1$; the latter follows from Theorem \ref{thm:gen} with $\Lambda = \ZZ_{\geq 0} \x \ZZ_{\geq 0}$. Hegarty \cite{Hegarty} showed that, for $d$ even, the limit $\rho^{s, d}$ is positive provided that it exists. However, it was previous unknown whether any of these limits exists.

Our next result provides some insight into the structure of a random subset $S \subseteq \set{0, 1, \dots, n}$ conditioned on $\lambda(S) \in \Lambda$. We argue that, except for the fringe elements of $S$ (i.e., the numbers close to $0$ or $n$), the middle elements are nearly unrestricted and independent from the fringe choices. The precise statement is found in Theorem~\ref{thm:middle}. This intuition was key to Martin and O'Bryant's proof \cite{MO} that $\rho_n$ is bounded below. It was also used by Miller, Orosz, and Scherinerman \cite{MOS} to construct a family of MSTD sets. However, previous work only applied the intuition to a relatively small proportion of all MSTD subsets. There has been no descriptions on what ``most'' MSTD sets look like. Our result is the first rigorous formulation of this common intuition. The techniques used in this paper have also inspired a new approach to a different problem on counting numerical semigroups of a given genus~\cite{Zhao:semigroup}.

For a uniformly random subset $S \subseteq\{0, 1, \dots, n\}$ conditioned on $\lambda(S) \in \Lambda$, our results imply that the middle segment of $S$ is close to being unrestricted. For instance, the probability that any ``middle'' element is in $S$ approaches $1/2$ as $n \to \infty$, thereby confirming (and generalizing) a conjecture of Miller, Orosz, and Scheinerman \cite{MOS}. Also, the expectation and variance of the size of $S$ are asymptotically the same as that of the binomial distribution on $n+1$ elements. The size distribution of $S$ also follow a central limit theorem.

This paper is organized as follows. We start by focusing exclusively on the MSTD problem. In Section \ref{sec:proof-limit} we show that the limit $\rho$ in Theorem \ref{thm:main} exists. In Section \ref{sec:computation} we elaborate on issues pertaining to computing lower and upper bounds for $\rho$. Next we move to the general case of subsets $S$ satisfying $\lambda(S) \in \Lambda$. In Section \ref{sec:extensions} we discuss how our methods for MSTD sets can be modified to prove Theorem \ref{thm:gen}. In Section \ref{sec:middle} we study the structure of a random set $S$ satisfying $\lambda(S) \in \Lambda$. Finally, in Section \ref{sec:conclusion} we offer some concluding remarks.

%%%%%%%%%%%%%%%%%%%%%%%%%%%%%%%%%%%%%%%%%%%%%%%%%%%%%%%%%%%%%%%%%%%%%%%%%%%%%%%%%%%%

\section{The limiting proportion of MSTD sets} \label{sec:proof-limit}

In this section, we show that proportion $\rho_n$ of MSTD sets converges to a limit. Although the proof contains a lot of the ingredients used in computing the limit, we defer to Section \ref{sec:computation} any details that are only relevant to the computation.

Let us give some intuition for our proof. Let $S$ be a ``typical'' subset of $\{0, 1, \dots, n\}$. As observed by Martin and O'Bryant \cite{MO}, except for elements near the ``fringe,'' most elements of $\{0, 1, 2, \dots, 2n\}$ can be represented as a sum of two elements of $S$ in a large number of ways. Consequently, these elements will ``typically'' be in the sum set. As Martin and O'Bryant put it, ``if we choose the `fringe' of $S$ cleverly, the middle of $S$ will be become largely irrelevant.'' 

The authors then proceed by manually fixing a particular choice of fringe for $S$, and thereby obtaining their lower bound for $\rho_n$. Unfortunately, fringe-fixing leads to very suboptimal lower bounds, since ``most'' MSTD sets do not have a particular fixed fringe profile.

Our idea is to let the fringe vary. For each particular fringe profile, we compute the proportion of subsets $S$ with the given fringe profile and the additional property that all the middle sums, namely those that are not completely controlled by the fringe, are in $S+S$. Then we can obtain the total proportion of MSTD subsets by summing over all candidate fringe profiles. Doing this leaves out those potential MSTD sets with some missing middle sum. Fortunately, as we will show, sets missing middle sums occupy a very small proportion of all subsets.

We begin by restricting ourselves to subsets $S \subseteq \{0, 1, \dots, n\}$ with $0, n \in S$, and then relax this constraint in Section \ref{sec:formula}.

\subsection{MSTD fringe pairs} \label{sec:fringe}

From now on, we use $[a,b]$ to denote the set $\{a, a+1, \dots, b\}$ if $a \leq b$, or the empty set otherwise.

Let $S \subseteq [0,n]$. When searching for fringe profiles candidates for $S$, we want the fringe alone to already generate more sums than differences. More precisely, suppose we fix $S \cap [0,k] = A$ and $(n - S) \cap [0, k] = B$. Then $(S+S) \cap [0,k]$ is completely controlled by $A$ and $(S + S) \cap [2n-k,2n]$ is completely controlled by $B$. Similarly, $(S - S) \cap (\pm [n-k, n])$ is completely controlled by $A$ and $B$. Suppose that we can choose the middle segment of $S$, i.e., $S\cap [k+1, n-k-1]$, so that every element of $[k+1, 2n-k-1]$ appears in $S + S$, then it would follow that $S$ is MSTD. So we would like to look for fringe profiles $(A,B)$ with the above properties. This is formalized in the following set of definitions. See Figure~\ref{fig:fringe} for a visual illustration.

\begin{figure}[ht]
  \centering
  \begin{tikzpicture}[xscale=.6, yscale=.4]
    \path[use as bounding box] (-13,-1) rectangle (13,10);
    \begin{scope}[yshift = 8cm]
      \node[above] at (-4.8, 1) {\scriptsize $0$};
      \node[above] at (4.9, 1) {\scriptsize $n$};
      \node[above] at (-3.2, 1) {\scriptsize $k$};
      \node[above] at (3.2, 1) {\scriptsize $n-k$};
      \node[below] at (-4, 0) {\scriptsize $A$};
      \node[below] at (4,0) {\scriptsize $n-B$};
      \draw (-5,0) rectangle (5,1);
      \draw[fill=gray] (-5,0) rectangle (-3,1);
      \draw[fill=gray] (5,0) rectangle (3,1);
      \node[left] at (-6, .5) {$S$};
    \end{scope}
    \begin{scope}[yshift = 4cm]
      \node[above] at (-9.8, 1) {\scriptsize $0$};
      \node[above] at (9.8, 1) {\scriptsize $2n$};
      \node[above] at (-8.2, 1) {\scriptsize $k$};
      \node[above] at (8.2, 1) {\scriptsize $2n-k$};
      \node[below] at (-9, 0) {\scriptsize $(A+A) \cap [0,k]$};
      \node[below] at (9,0) {\scriptsize $2n-(B+B)\cap [0,k]$};
      \draw (-10,0) rectangle (10,1);
      \draw[fill=gray] (-10,0) rectangle (-8,1);
      \draw[fill=gray] (10,0) rectangle (8,1);		
      \node[left] at (-11, .5) {$S+S$};
    \end{scope}
    \begin{scope}
      \node[above] at (-9.8, 1) {\scriptsize $-n$};
      \node[above] at (9.8, 1) {\scriptsize $n$};
      \node[above] at (-8.2, 1) {\scriptsize $-n+k$};
      \node[above] at (8.2, 1) {\scriptsize $n-k$};
      \node[below] at (-9, 0) {\scriptsize $-n + (A+B) \cap [0,k]$};
      \node[below] at (9,0) {\scriptsize $n - (A+B)\cap [0,k]$};
      \draw (-10,0) rectangle (10,1);
      \draw[fill=gray] (-10,0) rectangle (-8,1);
      \draw[fill=gray] (10,0) rectangle (8,1);	
      \node[left] at (-11, .5) {$S-S$};
    \end{scope}	
  \end{tikzpicture}

  \caption{The shaded areas are regions in $S$, $S+S$, and $S-S$ are completely controlled by the fringe $(A,B;k)$ of $S$.\label{fig:fringe}}
\end{figure}

\begin{definition} \label{def:fringe-pair}
A \emph{MSTD fringe pair} of order $k$ is a pair $(A, B)$ (also denoted $(A, B; k)$ to indicate the order), where $A$ and $B$ are both subsets of $[0,k]$, with $0 \in A$ and $0 \in B$, and satisfying
\[
	\abs{ (A+A) \cap [0,k] } + \abs{ (B+B) \cap [0,k] } > 2 \abs{ (A+B) \cap [0,k] }.
\]
\end{definition}

In Section \ref{sec:extensions} we consider a variation of fringe pairs to deal with generalizations of MSTD sets.

We impose the following partial order on the set of all MSTD fringe pairs: $(A, B; k) > (A', B'; k')$ if $k > k'$ and
\begin{equation}
  \label{eq:partial-order}
	A' = A \cap [0,k'], \quad B' = B \cap [0,k'], \quad	 [k'+1, k] \subseteq A+A, \quad [k'+1, k] \subseteq B + B.  
\end{equation}

\begin{definition}
A \emph{minimal MSTD fringe pair} is a MSTD fringe pair $(A, B; k)$ for which there does not exist another MSTD fringe pair $(A', B'; k')$ with $(A, B; k) > (A', B'; k')$.
\end{definition}

It is not hard to show that, to determine whether an MSTD fringe pair is minimal, it suffices to check \eqref{eq:partial-order} for $k'=k-1$. We use this fact in the computer search for minimal MSTD fringe pairs.

% It turns out that to determine whether an MSTD fringe pair is minimal, it suffices to check $k' = k-1$.

% \begin{lemma}
%   Let $(A,B;k)$ be an MSTD fringe pair. Then it is not minimal if and only if $k \in A+A, k \in B+B$ and $(A'',B'',k-1)$ is a fringe pair, where $A'' = A \cap[0,k-1]$ and $B'' = \cap [0, k-1]$.
% \end{lemma}

% \begin{proof}
%   The ``if'' direction is clear since we would have $(A, B; k) > (A', B'; k-1)$. For the ``only if'' direction, suppose that $(A, B; k)$ is not minimal, then there is some $k' < k$ such that \eqref{eq:partial-order} holds. If $k' = k-1$ then we are done. Otherwise, 
% \end{proof}

\begin{example}
  There are no MSTD fringe pairs of order less than 6. The minimal MSTD fringe pairs of order 6 are
  \[
  \begin{array}{ccc}
    A & B & k \\ \hline 
    \{0\} & \{0, 1, 3\} & 6 \\ 
    \{0\} & \{0, 2, 3\} & 6 \\ 
    \{0, 1, 3\} & \{0, 1, 2, 4\} & 6 \\ 
    \{0, 2, 3\} & \{0, 1, 2, 5\} & 6
  \end{array}
  \]
  as well as the four others obtained by switching $A$ and $B$. The minimal MSTD fringe pairs of order 7 are
  \[
  \begin{array}{ccc}
    A & B & k\\ \hline 
    \{0\} & \{0, 1, 3\} & 7 \\ 
    \{0\} & \{0, 2, 3\}  & 7 \\ 
    \{0\} & \{0, 1, 3, 4\} & 7 \\ 
    \{0\} & \{0, 1, 2, 5\} & 7 \\  
    \{0, 1, 3, 4\} & \{0, 1, 2, 5\} & 7
  \end{array}
  \]
  as well as the five others obtained by switching $A$ and $B$. There are ten non-minimal MSTD fringe pairs of order 7. They are
  \[
  \begin{array}{ccc}
    A & B & k \\ \hline 
    \{0,1,2,5\} & \{0,2,3,7\} & 7 \\ 
    \{0,7\} & \{0,1,3,7\} & 7 \\ 
    \{0,7\} & \{0,2,3,7\} & 7 \\  
    \{0,1,3,7\} & \{0,1,2,4,7\} & 7 \\ 
    \{0,2,3,7\} & \{0,1,2,5,7\} & 7
  \end{array}
  \]
  as well as the five others obtained by switching $A$ and $B$.
\end{example}

% Here we state some easy facts about minimal MSTD fringe pairs. Some of them are used in the computer search for such pairs. The proofs are mostly straight-forward, so we omit them.
% \begin{itemize}
% \item A minimal fringe pair $(A,B; k)$ cannot have $k \in A$ and $k \in B$ at the same time, but it is possible to have $k \in A$ and $k \notin B$.
% \end{itemize}

\begin{definition}
Let $S \subseteq [0,n]$. We say that $S$ is a \emph{rich MSTD set} with MSTD fringe pair $(A, B; k)$ if
\[
	2k < n, \quad
	S \cap [0,k] = A, \quad
	(n - S) \cap [0, k] = B, \quad\text{and}\quad
	[k+1, 2n-k-1] \subseteq S + S.
\]
The \emph{order} of the rich MSTD set $S$ is the smallest possible value of $k$ for which there exists such an MSTD fringe pair $(A, B; k)$.
\end{definition}

As expected, rich MSTD sets are MSTD, as we shall prove in a moment. We choose the name \emph{rich} because $S$ is rich in sums in the middle. Also, as we will see, they represent a rich collection of MSTD sets.

Next we prove some simple facts about rich MSTD sets and its MSTD fringe pairs. The goal is to show that we can count rich MSTD sets by going through the list of minimal MSTD fringe pairs. The proofs are mostly straightforward and they can be skipped if desired.

\begin{lemma} \label{lem:rich-MSTD}
A rich MSTD set is an MSTD set.
\end{lemma}

\begin{proof}
Let $S \subseteq [0, n]$ be a rich MSTD set with MSTD fringe pair $(A, B; k)$. We need to show that $\abs{S + S} > \abs{S - S}$. It suffices to show that
\begin{align}
	\label{eq:rich1} \abs{(S + S) \cap ([0,k] \cup [2n-k, 2n] )} &> \abs{(S - S) \cap ([-n, -n+k] \cup [n-k, n] )},   \\
	\label{eq:rich2} \text{and } \qquad \abs{(S + S) \cap [k+1, 2n-k-1]} &\geq \abs{(S - S) \cap [-n+k+1, n-k-1]}.
\end{align}
The inequality \eqref{eq:rich2} immediately follows from the requirement $[k+1, 2n-k-1] \subseteq S + S$. For \eqref{eq:rich1}, we note that
\begin{align*}
	(S + S) &\cap [0,k] = (A + A) \cap [0,k],
\\	(S + S) &\cap [2n-k, 2n] = ((n - B) + (n - B)) \cap [2n-k,2n] %= (2n - (B + B) )\cap [2n-k,2n] 
							= 2n - (B+B)\cap [0,k],
\\	(S - S) &\cap [-n, -n+k] = (A - (n-B))  \cap [-n, -n+k] = (A + B) \cap [0, k] - n,
\\	(S - S) &\cap [n-k, n] = ((n - B) - A) \cap [n-k, n] %= (n - (A + B)) \cap [n-k, n]  
							= n - (A + B) \cap [0,k].
\end{align*}
And hence the sizes of the above four sets are $\abs{(A + A) \cap [0,k]}, \abs{(B+B)\cap [0,k]}, \abs{(A + B) \cap [0,k]}$, and $\abs{(A + B) \cap [0,k]}$, respectively. Then \eqref{eq:rich1} follows from $(A, B; k)$ being an MSTD fringe pair.
\end{proof}

A rich MSTD set may have many choices for its fringe pair. The following lemma shows that the set of MSTD fringe pairs of a particular rich MSTD set forms a chain in the partial order.

\begin{lemma} \label{lem:fringe-order1}
Let $S \subseteq [0, n]$ be a rich MSTD set. Let $(A, B; k)$ and $(A', B'; k')$ be two MSTD fringe pairs of $S$. If $k = k'$, then $(A, B; k) = (A', B'; k')$. If $k > k'$, then $(A, B; k) > (A', B'; k')$.
\end{lemma}

\begin{proof}
If $k = k'$, then $A = A' = S \cap [0, k]$ and $B = B' = (n - S) \cap [0, k]$. So $(A, B; k) = (A', B'; k')$.

If $k > k'$, then $A' = S \cap [0, k'] = A \cap [0, k']$, $B' = (n - S) \cap [0, k'] = B \cap [0, k']$. Since $S$ is rich with fringe pair $(A', B', k')$, we see that $[k'+1, 2n - k'-1] \subseteq S + S$. The sum in $[0, k]$ can only come from a sum of two elements in $[0, k]$, so that $[k' + 1, k] \subseteq A+A$. Similarly, $[k'+1, k] \subseteq B+B$. Hence $(A, B; k) > (A', B'; k')$.
\end{proof}

Thus, for a rich MSTD set of order $k$, we can speak of its \emph{minimal} MSTD fringe pair, which necessarily has order $k$.

\begin{lemma} \label{lem:fringe-order2}
Let $S \subseteq [0, n]$ be a rich MSTD set. Let $(A, B; k)$ be the minimal MSTD fringe pair of a rich MSTD set $S$. Then $(A, B; k)$ is minimal in the partial ordering of all MSTD fringe pairs. Also, for every $k < k' < n/2$, $(A', B'; k')$ is also an MSTD fringe pair of $S$, where $A' = S \cap [0, k]$ and $B' = (n - S) \cap [0, k]$, and every MSTD fringe pairs of $S$ have this form.
\end{lemma}

\begin{proof}
Suppose that $(A, B; k)$ is not a minimal MSTD fringe pair, so that we have $(A', B'; k') < (A, B; k)$. Then $A' = A \cap [0, k'] = S \cap [0, k']$ and $B' = B \cap [0, k'] = (n-S) \cap [0, k']$. Also, $[k'+1, k]$ is contained in $A+A$ and $B+B$, and $[k+1, 2n-k-1] \subseteq S+S$ (since $S$ is rich of order $k$), so that $[k'+1, 2n-k'-1] \subseteq S+S$. Hence $(A', B'; k')$ is also a fringe pair of $S$, thereby contradicting the the choice of $(A, B; k)$ as the minimal MSTD fringe pair of a rich MSTD set $S$.

%We need to check that
%\begin{equation} \label{eq:fringe-order2}
%	\abs{ (A'+A') \cap [0,k'] } + \abs{ (B'+B') \cap [0,k] } > 2 \abs{ (A'+B') \cap [0,k'] }.
%\end{equation}
For the second claim, where $k < k'$, we see that $[k+1, k']$ is contained in $A'+A'$ and $B'+B'$ as $[k+1, 2n-k-1] \subseteq S+S$. Since $(A, B; k)$ is an MSTD fringe pair, we have
\begin{align*}
	\abs{ (A'+A') \cap [0,k'] } + \abs{ (B'+B') \cap [0,k] }
&	= \abs{(A+A)\cap[0,k]} + \abs{(B+B)\cap[0,k]} + 2(k'-k)
\\&	> 2\abs{(A+B)\cap[0,k]} + 2(k'-k)
\\&	\geq 2 \abs{ (A'+B') \cap [0,k'] }
\end{align*}
Hence $(A', B'; k')$ is also an MSTD fringe pair. The rest of the lemma is clear.
\end{proof}

Therefore, we can count rich MSTD sets by their minimal MSTD fringe pairs.

%As a consequence, if $(A, B; k)$ is the minimal MSTD fringe pair of a rich MSTD set $S$, then for every $k < k' < n/2$, $(A', B'; k')$ is an MSTD fringe pair of $S$, where $A' = S \cap [0, k]$ and $B' = n - S \cap [n, n-k]$, and all the MSTD fringe pairs of $S$ have this form. We will count rich MSTD sets by their minimal MSTD fringe pairs. By only considering minimal MSTD fringe pairs, we do not run into the risk of overcounting rich MSTD sets.

\subsection{Semi-rich sets} \label{sec:semi-rich}

We are interested in counting the number of rich MSTD sets with a given MSTD fringe pair. It turns out that we can divide this problem into two halves: the front half and back half. In this section we show how to compute the relevant limiting proportions for each half. In the next section we show how to put the two halves together.

%We introduce the notion of a semi-rich set, which, very roughly speaking, is ``half'' of a rich set.

\begin{definition}
We say that $T \subseteq [0, n]$, where $n \geq k$, is a \emph{$k$-semi-rich} set if $[k+1, n] \subseteq T + T$. We say that $T$ has \emph{prefix} $(A;k)$ where $A = T \cap [0, k]$.
\end{definition}

For $n \geq k$ and $A \subseteq [0,k]$ (with $0 \in A$), let
	%originally, F(A; k, n) = \sigma_n(A; k) 2^n 
\begin{equation} \label{eq:sigma-def}
	\sigma_n(A; k) = 2^{-n} \abs{\set{ T \subseteq [0,n] : T \cap [0, k] = A, [k+1, n] \subseteq T + T }}.
\end{equation}
In other words, $\sigma_n(A;k)$ is the probability that a uniformly random subset $S \subseteq [0,n]$ (conditioned on $0 \in S$) is $k$-semi-rich with prefix $(A;k)$. In this section, we show that $\sigma_n(A; k)$ converges to a limit and give a formula for computing this limit.

\begin{proposition} \label{prop:sigma-exists}
For every $A \subseteq [0, k]$ with $0 \in A$, the limit
\[
	\sigma(A; k) = \lim_{n \to \infty} \sigma_n(A; k)
\]
exists and is positive.
\end{proposition}

\begin{proof}
We compute the size of the collection in \eqref{eq:sigma-def} by considering the complement. We know that
\begin{equation}\label{eq:F-compl}
	\sigma_n(A; k) 2^n = 2^{n-k} - \abs{\set{ T \subseteq [0,n] : T \cap [0, k] = A, [k+1, n] \not\subseteq T + T }}.
\end{equation}
Observe that the set on the RHS can be partitioned by the smallest element of $[k+1, n]$ not in $T + T$, that is,
\begin{multline*}
	\{ T \subseteq [0,n] : T \cap [0, k] = A, [k+1, n] \not\subseteq T + T \}
\\	= \biguplus_{j > k} \{ T \subseteq [0,n] : T \cap [0, k] = A, [k+1, j-1] \subseteq T + T, j \notin T + T \}
\end{multline*}
where $\uplus$ denotes disjoint union. We introduce the following quantity for $j > k$:
\[
	G_j(A; k) = \abs{\set{ T \subseteq [0,j] : T \cap [0, k] = A, [k+1, j-1] \subseteq T + T, j \notin T + T }}.
\]
Then, for $k < j \leq n$,
\begin{align*}
	&\abs{\set{ T \subseteq [0,n] : T \cap [0, k] = A, [k+1, j-1] \subseteq T + T, j \notin T + T }}
\\&	\quad = 2^{n-j} \cdot \abs{\set{ T \subseteq [0,j] : T \cap [0, k] = A, [k+1, j-1] \subseteq T + T, j \notin T + T }}
\\&	\quad = G_j(A; k) 2^{n-j},
\end{align*}
since $T \cap [j+1, n]$ can be chosen arbitrarily. It follows from \eqref{eq:F-compl} that
\[
	\sigma_n(A; k) 2^n = 2^{n-k} - \sum_{j = k+1}^n G_j(A; k) 2^{n-j}.
\]
So
\[
	\sigma_n(A; k)  = 2^{-k} - \sum_{j = k+1}^n G_j(A; k) 2^{-j},
\]
and hence
\begin{equation} \label{eq:sum}
	\sigma(A; k) = \lim_{n \to \infty} \sigma_n(A; k)
	= 2^{-k} - \sum_{j = k+1}^\infty G_j(A; k) 2^{-j}.
\end{equation}
In particular, the limit exists since the quantities $G_j(A; k)$ and $\sigma_n(A; k)$ are all non-negative. The argument for $\sigma(A; k) > 0$ is very similar to the arguments in $\cite{MO}$, so we only sketch the idea. Basically, if we choose a sufficiently large $\ell$ (depending on $k$) and require that $[k+1, \ell] \subseteq T$, and then choose $T \cap [\ell+1, n]$ randomly, then there is a positive lower bounded probability that $[k+1, n] \subseteq T + T$, thereby making $T$ semi-rich (the idea is very similar to Lemma~\ref{lem:missing-middle}).
\end{proof}

\subsection{Rich MSTD sets with a given MSTD fringe pair} \label{sec:rich}

Fix an MSTD fringe pair $(A,B;k)$. As $n \to \infty$, what proportion of the subsets of $[0,n]$ are rich MSTD sets with MSTD fringe pair $(A,B;k)$? In this section, we show that the answer is simply the product of the proportions of $k$-semi-rich sets with prefix $(A; k)$ and $(B; k)$ respectively. 

The intuition here is that, for large $n$ and a uniform random subset $S \subseteq [0, n]$, with very high probability every element in $[n/2, 3n/2]$ appears in the sum set $S+S$. So we are mostly concerned with ensuring that each half of $S$ is semi-rich.

For an MSTD fringe pair $(A,B;k)$, and an integer $n > 2k$, let
\begin{equation}
  \label{eq:rho-fringe-pair}
  \rho_n(A,B;k) = 2^{-n+1} \abs{\set{ S \subseteq [0,n] : S \cap [0,k] = A, (n-S) \cap[0,k] = B, 
      [k+1, 2n-k-1] \subseteq S + S}}.
\end{equation}
In other words, $\rho_n(A,B;k)$ is the probability that a uniformly chosen random subset $S \subseteq [0,n]$ (conditioned on $0,n \in S$) is a rich MSTD set with MSTD fringe pair $(A,B;k)$. The following proposition formalizes the above intuition.

\begin{proposition} \label{prop:product}
As $n \to \infty$, $\rho_n(A, B; k)$ approaches a limit $\rho(A, B; k)$, and
\[
	\rho(A, B; k) = \sigma(A; k) \sigma(B; k).
\]
\end{proposition}

\begin{proof}
In this proof, assume that $n$ is sufficiently large. Let $m =
\floor{n/2}$. If a subset $S \subseteq [0, n]$ is a rich MSTD subset with MSTD fringe pair $(A,B; k)$, then it follows that $S \cap [0, m]$ is a $k$-semi-rich subset of $[0,m]$ with prefix $(A; k)$, and $(n - S) \cap [0, n-m-1]$ is a $k$-semi-rich subset of $[0, n-m-1]$ with prefix $(B;k)$. Thus we have
\begin{equation} \label{eq:prop:product1}
	\rho_n(A, B; k) 2^{n-1} \leq \sigma_m(A; k) 2^m \cdot \sigma_{n-m-1}(B;k) 2^{n-m-1} = \sigma_m(A; k) \sigma_{n-m-1}(B;k) 2^{n-1}.
\end{equation}
The difference $\sigma_m(A; k) \sigma_{n-m-1} (B; k) 2^{n-1} - \rho_n(A, B; k) 2^{n-1}$ counts the collection of subsets of $[0, n]$ which, among other things, have the property that some element in $[m+1, n+m]$ is missing from $S + S$. It is easy to see that the number of subsets $S \subseteq [0, n]$ satisfying $j \notin S+S$ is precisely $3^{\floor{(j'+1)/2}}\cdot 2^{n-j'}$ where $j' = j$ if $0 \leq j \leq n$ and $j' = 2n - j$ if $n < j \leq 2n$. So, if $j \in [m+1, n+m]$, then the number of subsets $S \subseteq [0,n]$ with $j \notin S + S$ is at most $3^{m/2} 2^{n-m} \leq 3^{n/4}2^{n/2 + 1}$ (recall that $m = \floor{n/2}$). Therefore,
\begin{align} 
	\sigma_m(A; k) \sigma_{n-m-1}(B;k) 2^{n-1} - \rho_n(A, B; k) 2^{n-1} 
&	\leq	\abs{\set{ S \subseteq [0, n] : [m+1, n+m] \not\subseteq S + S }} \nonumber
\\&	\leq n \cdot 3^{n/4}2^{n/2 + 1}. \label{eq:prop:product2}
\end{align}
Combining \eqref{eq:prop:product1} and \eqref{eq:prop:product2} we obtain
\[
	\sigma_m(A; k) \sigma_{n-m-1}(B;k) - n \cdot 3^{n/4}2^{-n/2 + 2} \leq \rho_n(A, B; k) \leq \sigma_m(A; k) \sigma_{n-m-1}(B; k).
\]
Letting $n \to \infty$ gives
\[
	\lim_{n \to \infty} \rho_n(A, B; k) = \lim_{n \to \infty} \sigma_{\floor{n/2}}(A; k) \sigma_{n-\floor{n/2}-1}(B; k)
										= \sigma(A; k) \sigma(B; k),
\]
thereby proving the lemma.
%So
%\[
%	F(A; k, m)F(B; k, n-m-1) - n \cdot 3^{n/4}2^{n/2 + 1} \leq P(A, B; k, n) \leq F(A; k, m)F(B; k, n-m-1).
%\]
%Then
%\begin{align*}
%	&F(A; k, m)F(B; k, n-m-1)2^{-n+1} - n \cdot 3^{n/4}2^{-n/2+2}  \\
%	&\qquad \leq P(A, B; k, n) 2^{-n+1} \leq F(A; k, m)F(B; k, n-m-1)2^{-n+1}.
%\end{align*}
%Letting $n \to \infty$, we obtain that
%\begin{align*}
%	\rho(A, B; k) 
%&	= \lim_{n \to \infty} P(A, B; k, n) 2^{-n+1}
%\\&	= \lim_{n \to \infty} F(A; k, \floor{n/2})F(B; k, n-\floor{n/2}-1)2^{-n+1}
%\\&	= \sigma(A; k) \sigma(B; k),
%\end{align*}
%as desired.
\end{proof}
%
%Therefore, using the estimates $f_-(A; k)$ and $f_+(A; k)$ for $\sigma(A; k)$, we can obtain corresponding lower and upper estimates for $\rho(A, B; k)$:
%\[
%	\rho_-(A, B; k) = \frac 1 4 f_-(A; k)f_-(B; k), \quad \text{and} \quad
%	\rho_+(A, B; k) = \frac 1 4 f_+(A; k)f_+(B; k).
%\]

\subsection{Almost all MSTD sets are rich} \label{sec:most-rich}

Previously we considered the proportion of rich MSTD sets with a particular MSTD fringe pair. By summing over all minimal MSTD fringe pairs, we obtain the proportion of rich MSTD sets. In this section, we show that, in some sense, almost all MSTD sets are rich, so that the limiting proportion of MSTD sets equals the limiting proportion of rich MSTD sets.

The intuition, as before, is that there is a diminishingly small probability that any ``middle'' sum or difference is missing. We can quantify this observation through the following two lemmas.

\begin{lemma} \label{lem:missing}
Let $S$ be a uniform random subset of $[0, n]$ containing $0$ and $n$.
\begin{enumerate}
	\item[(a)] If $s \in [1, n-1]$, then
	\[
		\PP\left\{ s \notin S + S \right\}
		= \left\{\begin{array}{ll}
			\frac{1}{2}\(\frac 3 4\)^{(s-1)/2}, & \text{if } s \text{ is odd} \\
			\frac{1}{4} \( \frac 3 4\)^{(s-2)/2}, & \text{if } s \text{ is even}
			\end{array}
			\right\} \leq \frac{1}{2}\(\frac 3 4\)^{(s-1)/2}.
	\]
	And if $s \in [n+1, 2n-1]$, then $\PP\{s \notin S + S\} = \PP \{ 2n - s \notin S + S\}$.
	
	\item[(b)] If $d$ is an integer with $n/2 < d < n$, then
	\[
		\PP\left\{ d \notin S - S \right\} = \frac 1 4 \(\frac{3}{4}\)^{n-d-1}.
	\]
	If $0 < d \leq n/2$, then
	\[
		\PP\left\{ d \notin S - S \right\} \leq \(\frac{3}{4}\)^{(n-1)/3}.
	\]
	Finally, $\PP\{d \notin S - S \} = \PP\{-d \notin S - S\}$.
\end{enumerate}
\end{lemma}

We omit the easy proof of Lemma \ref{lem:missing} since very similar results can be found in \cite[Sec.~2]{MO}. We also used similar ideas in the proof of Proposition \ref{prop:product}.

\begin{lemma} \label{lem:missing-middle}
Let $n$ and $\kbar$ be positive integers with $n > 2\kbar$. Let $S$ be a uniform random subset of $[0, n]$ containing $0$ and $n$. Then
\[
	\PP\left\{ [\kbar+1, 2n-\kbar-1] \not\subseteq S + S \right\}
		\leq \frac{(3/4)^{\kbar/2}}{1 - \frac{\sqrt 3}2},
\]
and
\[
	\PP\left\{ [-n+\kbar+1, n-\kbar-1] \not\subseteq S - S \right\} \leq 2 \(\frac 3 4\)^{\kbar} + (n+1) \(\frac 3 4\)^{(n-1)/3}.
\]
\end{lemma}

\begin{proof}
In each case, apply the union bound, use Lemma \ref{lem:missing}, and then sum a geometric series.
\end{proof}

We also state a variation Lemma \ref{lem:missing-middle} where we drop the restriction that $S$ contains $0$ and $n$. The proof is very similar so we omit it. This lemma will be used in later sections.

\begin{lemma} \label{lem:not-affluent}
Let $n$ and $\kbar$ be positive integers with $n > 2\kbar$. Let $S$ be a uniform random subset of $[0, n]$. Then
\[
	\PP\left\{ [\kbar+1, 2n-\kbar-1] \not\subseteq S + S \right\}
		\leq \frac{3(3/4)^{\kbar/2}}{2 - \sqrt 3}
\]
and
\[
	\PP\left\{ [-n+\kbar+1, n-\kbar-1] \not\subseteq S - S \right\} 
		\leq 8 \( \frac{3}{4}\)^{\kbar+2} + (n+1) \(\frac{3}{4}\)^{(n-1)/3}.
\]
\end{lemma}

The take-away point from the above two lemmas is that by forcing $\kbar$ to be large, we can make the probability that any middle sum or difference is missing to be negligible. In other words,
\begin{align*}
	\lim_{\kbar \to \infty} \limsup_{n \to \infty} \PP\left\{ [\kbar+1, 2n-\kbar-1] \not\subseteq S + S \right\} &= 0, \\
	\lim_{\kbar \to \infty} \limsup_{n \to \infty} \PP\left\{ [-n+\kbar+1, n-\kbar-1] \not\subseteq S - S \right\} &= 0.
\end{align*}

Now we state the result that formalizes the statement that ``almost all MSTD sets are rich.'' For now, we restrict ourselves to MSTD sets $S \subseteq [0, n]$ containing $0$ and $n$. Let
\[
	\rho_{*n} = 2^{-n+1} \abs{\set{S \subseteq [0, n] : 0, n \in S, \text{ and } S \text{ is MSTD} }}.
\]
We put the asterisk in the subscript to indicate that $0, n \in S$ because we need to reserve the superscript space for later.

\begin{proposition} \label{prop:almost-all}
As $n \to \infty$, $\rho_{*n}$ converges to a limit $\rho_*$, and
\[
	\rho_* = \sum_{(A, B; k)} \rho(A, B; k)
\] 
where the sum is taken over all minimal MSTD fringe pairs $(A, B; k)$.
\end{proposition}

\begin{proof}
Fix $\kbar$ a positive integer. We start by considering only MSTD fringe pairs of order at most $\kbar$. In the last step of the proof we let $\kbar \to \infty$.

Assume that $n$ is sufficiently large. If $S$ is a uniform random subset of $[0, n]$ containing $0$ and $n$, then $\rho_{*n}$ is the probability that $S$ is MSTD. Since rich MSTD sets of order at most $\kbar$ form a subset of all MSTD sets, we have
\begin{equation} \label{eq:lower-approx}
	\sum_{\substack{(A, B; k) \\ k \leq \kbar}} \rho_n(A, B; k) \leq \rho_{*n}.
\end{equation}
Unless otherwise specified, such sums are always assumed to be taken over minimal MSTD fringe pairs. Note that the sum has finitely many terms.

Let $S \subseteq [0, n]$ be an MSTD set containing $0$ and $n$. Let $A = S \cap [0, \kbar]$ and $B = (n- S) \cap [0, \kbar]$ be the fringe sets as usual. Suppose that $S$ is not a rich MSTD set of order at most $\kbar$ (meaning that either $S$ is not rich, or $S$ is rich with order greater than $\kbar$). There are two possibilities
\begin{description}
	\item[Case 1.] $(A,B; \kbar)$ is not an MSTD fringe pair. Then 
				\[
					\abs{(A+A) \cap [0, \kbar]} + \abs{(B+B) \cap [0, \kbar]} \leq \abs{(A+B) \cap [0,\kbar]}.
				\]					
					Since $S$ is an MSTD set, $S-S$ must be missing some difference in $[-n+\kbar+1, n-\kbar+1]$ (c.f.~proof of Lemma \ref{lem:rich-MSTD}).
	\item[Case 2.] $(A, B; \kbar)$ is an MSTD fringe pair, but $S \subseteq [0, n]$ is not a rich MSTD set of $\kbar$, i.e., 
	$S+S$ is missing some sum in $[\kbar+1, 2n-\kbar-1]$.
	%$\abs{(A+A) \cap [0, \kbar]} + \abs{(B+B) \cap [0, \kbar]} > \abs{(A+B) \cap [0,\kbar]}$					
\end{description}

In both cases, $S$ is missing a middle sum or a middle difference. By Lemma \ref{lem:missing-middle}, we have
\begin{align*}
	0 &\leq \rho_{*n} -\sum_{\substack{(A, B; k) \\ k \leq \kbar}} \rho_n(A, B; k)
	  \\&\leq \PP\left\{ [\kbar+1, 2n-\kbar-1] \not\subseteq S + S \right\} + \PP\left\{ [-n+\kbar+1, n-\kbar-1] \not\subseteq S - S \right\}.
	  \\&\leq \frac{(3/4)^{\kbar/2}}{1 - \frac{\sqrt 3}2} + 2 \(\frac 3 4\)^{\kbar} + (n+1) \(\frac 3 4\)^{(n-1)/3}.
\end{align*}
Let $n \to \infty$ and we get
\begin{equation} \label{eq:almost-all-bounds}
	\sum_{\substack{(A, B; k) \\ k \leq \kbar}} \rho(A, B; k)
	\leq \liminf_{n \to \infty} \rho_{*n}
	\leq \limsup_{n \to \infty} \rho_{*n}
	\leq \(\sum_{\substack{(A, B; k) \\ k \leq \kbar}} \rho(A, B; k) \) +  \frac{(3/4)^{\kbar/2}}{1 - \frac{\sqrt 3}2} + 2 \(\frac 3 4\)^{\kbar}.
\end{equation}
Let $\kbar \to \infty$ and we get
\[
	\rho_{*} = \lim_{n \to \infty} \rho_{*n} =  \sum_{(A, B; k)} \rho(A, B; k). \qedhere
\]
\end{proof}

\subsection{The proportion of MSTD sets} \label{sec:formula}

In this section we remove the restriction that $0, n \in S$. Recall that $\rho_n$ is the probability that a uniform random subset of $[0,n]$ is an MSTD set.

\begin{lemma} \label{lem:p*=p}
$\displaystyle \lim_{n \to \infty} \rho_n = \rho_{*}$.
\end{lemma}

\begin{proof}
Fix $\epsilon > 0$. Choose an $N$ so that $\abs{\rho_{*m} - \rho_*} < \epsilon$ for all $m > N/3$. Let $S$ be a random subset of $[0, n]$, where $n > N$. Let $E$ denote the event that $\min S < n/3$ and $\max S > 2n/3$. So $\PP(E) = (1 - 2^{-\floor{n/3}+1})^2$. If $E$ occurs, then the probability that $S$ is MSTD is $\epsilon$-close to $p_*$. It follows that
\[
	\(1 - 2^{-\floor{n/3}+1}\)^2 (\rho_* - \epsilon) < \rho_n < \(1 - 2^{-\floor{n/3}+1}\)^2 (\rho_* + \epsilon) + 1 - \(1 - 2^{-\floor{n/3}+1}\)^2.
\]
for $n > N$. Let $n \to \infty$ and we get
\[
	\rho_* - \epsilon  \leq \liminf_{n \to \infty} \rho_n \leq \limsup_{n \to \infty} \rho_n \leq \rho_* + \epsilon.
\]
Since $\epsilon$ can be made arbitrarily small, we have
\[
	\lim_{n \to \infty} \rho_n = \rho_*. \qedhere
\]
\end{proof}

Combining Propositions \ref{prop:product}, \ref{prop:almost-all} and Lemma~\ref{lem:p*=p}, we obtain the following formula for the density of MSTD sets.

\begin{proposition} \label{prop:formula}
The density of MSTD sets satisfy
\[
	\rho = \lim_{n\to\infty} \rho_n = \sum_{(A, B; k)} \rho(A, B; k) = \sum_{(A, B; k)}\sigma(A,k) \sigma(B; k)
\]
where the sum is taken over all minimal MSTD fringe pairs $(A, B; k)$.
\end{proposition}

In particular, we have proven the existence of the limit in Theorem \ref{thm:main}. Proposition \ref{prop:formula} also gives the formula that we will use to compute $\rho$.

\section{Computing the limit} \label{sec:computation}

In this section we explain how to compute lower and upper bounds for $\rho$. Our method could, in principle, be used to derive bounds of arbitrary precision, although in practice the computation time increases exponentially with desired precision. We start with a description of the method to compute the estimate to $\rho$. Our numerical results can be found at the end of this section.

Our computation consists of the following steps. The functions $\sigma$ and $\rho$ were defined in Sections~\ref{sec:semi-rich} and \ref{sec:rich}, respectively.

\begin{enumerate}
\item Fix a $\kbar$. Find all minimal MSTD fringe pairs of order up to $\kbar$.
\item For each $(A, B; k)$ found in step 1, compute lower and upper bounds for $\sigma(A; k)$ and $\sigma(B; k)$.
\item Add up the lower and upper bounds for $\rho(A, B; k) = \sigma(A; k)\sigma(B; k)$ for all $(A, B; k)$ found in step 1.
\end{enumerate}

The variables $\kbar$, $\jbar$, and $h_k$ are all computational parameters, viewed as inputs to the computation. Each variable represents the  extent of some complete search. In general, larger values of these parameters give better numerical results but also increases running time.

\subsection{Generating minimal MSTD fringe pairs}

All the minimal MSTD fringe pairs of order $k$ can generated by a complete search through all pairs subsets of $[0, k]$, for each $k$ up to $\kbar$. That is, we generate a list of all pairs of subsets $A, B \subseteq [0,k]$ satisfying
\begin{itemize}
\item $0 \in A, 0 \in B$;
\item $\abs{ (A+A) \cap [0,k] } + \abs{ (B+B) \cap [0,k] } > 2 \abs{ (A+B) \cap [0,k] }$;
\item The following statements are not all true: $k \in A + A$, $k \in B + B$, $\abs{ (A+A) \cap [0,k-1] } + \abs{ (B+B) \cap [0,k-1] } > 2 \abs{ (A+B) \cap [0,k-1] }$.
\end{itemize}
The first two items correspond to $(A, B; k)$ being an MSTD fringe pair, while the third item corresponds to minimality.

\subsection{Estimating $\sigma(A; k)$} \label{sec:estimate-sigma}

Recall that $\sigma(A; k)$ is the density of semi-rich sets with prefix $(A; k)$. The methods used here to compute lower and upper bounds to $\sigma(A; k)$ build on the results developed earlier in Section \ref{sec:semi-rich}.

The key formula is \eqref{eq:sum}, which we reproduce here for convenience:
\begin{equation} \label{eq:sum2}
	\sigma(A; k) = 2^{-k} - \sum_{j = k+1}^\infty G_j(A; k) 2^{-j}
\end{equation}
where
\[
	G_j(A; k) = \abs{\set{ T \subseteq [0,j] : T \cap [0, k] = A, [k+1, j-1] \subseteq T + T, j \notin T + T }}.
\]
The computation consists of the following steps. Here $\jbar$ is a computational parameter.
\begin{enumerate}
\item Compute the terms $G_j(A; k)$ in \eqref{eq:sum2} for all $j$ satisfying $k < j \leq \jbar$ to obtain an upper bound to $\sigma(A;k)$ by using a partial sum.
\item Upper bound the trailing sum $\sum_{j > \jbar} G_j(A; k) 2^{-j}$ in \eqref{eq:sum2} to obtain a lower bound to $\sigma(A;k)$.
\end{enumerate}
In this section, we describe how to produce two numbers $\sigma_-(A;k)$ and $\sigma_+(A;k)$ such that
\[
\sigma_-(A;k) \leq \sigma(A;k) \leq \sigma_+(A;k).
\]
By increasing our computational parameters, we could, in principle, make the two estimates $\sigma_-(A;k)$ and $\sigma_+(A;k)$ arbitrarily close to the true value $\sigma(A;k)$. Unfortunately, the cost of computation increases prohibitly with desired precision level.

\subsubsection{Upper estimate of $\sigma(A;k)$}

Each individual term $G_j(A; k)$ can be computed by a complete search. For each minimal MSTD fringe pair $(A;k)$, we shall compute $G_j(A;k)$ for all $j$ satisfying $k < j \leq \jbar$. Our upper bound to $\sigma(A;k)$ is then given by
\begin{equation}
  \label{eq:sigma-upper}
\sigma_+(A;k) = 2^{-k} - \sum_{j = k+1}^{\jbar} G_j(A; k) 2^{-j}  
\end{equation}

\subsubsection{Lower estimate of $\sigma(A;k)$}

To determine a lower estimate of $\sigma(A;k)$, we need an effective upper bound for the following the trailing terms in \eqref{eq:sum2}:
\begin{equation} \label{eq:trailing}
	\sum_{j > \jbar} G_j(A; k) 2^{-j}.
\end{equation}
In computing an upper bound to \eqref{eq:trailing}, we do not explicit compute the exact values of any additional $G_j(A;k)$ terms. We obtain an upper bound through the following series of lemmas:

\begin{lemma} \label{lem:error1}
Let $A \subseteq [0, k]$. If $2k < j$, then
\begin{equation} \label{eq:error1a}
	G_j(A; k) \leq 2^{k+1 - \abs{A}} \cdot 3^{\floor{(j-2k-1)/2}}
\end{equation}
and if $k < j \leq 2k$, then $G_j(A; k) = 0$ if $j \in A + A$, and otherwise
\begin{equation} \label{eq:error1b}
	G_j(A; k) \leq 2^{j - k - \abs{A \cap [0, j-k-1]}}.
\end{equation}
\end{lemma}

\begin{proof}
In both cases, the bound simply uses the fact that
\begin{equation} \label{eq:error1c}
	G_j(A; k) \leq \abs{\set{ T \subseteq [0,j] : T \cap [0, k] = A, j \notin T + T }}.
\end{equation}
It can be easily checked that the RHS of the \eqref{eq:error1c} equals to the RHS expression in \eqref{eq:error1a} and \eqref{eq:error1b} in the respective cases.
\end{proof}

\begin{lemma} \label{lem:error2}
Let $A \subseteq [0, k]$ and $2k < \ell$. Then
\[
	\sum_{j = \ell}^\infty G_j(A; k) 2^{-j} \leq 
	\begin{cases}
	2^{k+2 - \abs{A} - \ell} \cdot 3^{-k + \frac {\ell+1} 2} & \text{if } \ell \text{ is odd,} \\
	5 \cdot 2^{k+2 - \abs{A} - \ell} \cdot 3^{-k-1 + \frac {\ell} 2} & \text{if } \ell \text{ is even.}
	\end{cases}	
\]
\end{lemma}

\begin{proof}
This follows from applying Lemma~\ref{lem:error1} to each term in the infinite sum, and then summing a geometric series:
\[
	\sum_{j = \ell}^\infty G_j(A; k) 2^{-j}
	\leq \sum_{j = \ell}^\infty  2^{k+1 - \abs{A}-j} \cdot 3^{\floor{(j-2k-1)/2}}
	= \frac{2^{k+1 - \abs{A}-\ell} \cdot 3^{\floor{(\ell-2k-1)/2}} + 2^{k - \abs{A}-\ell} \cdot 3^{\floor{(\ell-2k)/2}}}{1 - \frac{3}{4}}.
\]
The last expression above equals to the upper bound given in the lemma.
\end{proof}

Lemma \ref{lem:error2} is sufficient in providing an upper bound to \eqref{eq:sum}. However, the bound turns out to be somewhat weak. That is, in theory we already have the tools to evaluate the limit in \eqref{eq:sum} to arbitrary precision, but we would like an more efficient way of upper bounding the trailing error terms \eqref{eq:trailing}. This issue is handled by the following lemma.

\begin{lemma} \label{lem:error3}
Let $k < h < j$ and $A \subseteq [0, k]$. Let $\c B_{h}(A;k)$ denote the set of all $B \subseteq [0,h]$ satisfying $B \cap [0, k] = A$ and $[k+1, h] \subseteq B + B$. Then
\[
	G_j(A; k) = \sum_{B \in \c B_h(A;k)} G_j(B, h).
\]
\end{lemma}

\begin{proof}
The lemma follows from taking the cardinality of
\begin{multline*}
  \{ T \subseteq [0,j] : T \cap [0, k] = A, [k+1, j-1] \subseteq T + T, j \notin T + T \} \\
  = \biguplus_{B \in \c B_h(A;k)} \{ T \subseteq [0,j] : T \cap [0, h] = B, [h+1, j-1] \subseteq T + T, j \notin T + T \}. \qedhere
\end{multline*}
\end{proof}

We will use Lemma~\ref{lem:error3} in way that allows $h$ to vary with $k$. Let $h_k$ be a computational parameter, one for each $k$.

Our method of computing the upper bound to \eqref{eq:trailing} combines Lemmas \ref{lem:error1}, \ref{lem:error2}, and \ref{lem:error3}. In other words, let $\ol G_j(A;k)$ denote the upper bound to $G_j(A;k)$ given in Lemma~\ref{lem:error1}, and denote the upper bound in Lemma~\ref{lem:error2} by
\[
GT_\ell(A;k) = 	\begin{cases}
	2^{k+2 - \abs{A} - \ell} \cdot 3^{-k + \frac {\ell+1} 2} & \text{if } \ell \text{ is odd,} \\
	5 \cdot 2^{k+2 - \abs{A} - \ell} \cdot 3^{-k-1 + \frac {\ell} 2} & \text{if } \ell \text{ is even.}
	\end{cases}	
\]
Then we have
\[
\sum_{j > \jbar} G_j(B; h) 2^{-j}
\leq GT_{\max\set{2h+1, \jbar+1}}(A;k) + \sum_{\jbar < j \leq 2h} \ol G_j(B; h) 2^{-j}.
\]
Then
\begin{align} 
	\sum_{j > \jbar} G_j(A; k) 2^{-j}
	&= \sum_{B \in \c B_{h_k}(A;k)} \sum_{j > \jbar} G_j(B, h_k)  2^{-j} \nonumber
        \\	&\leq  \sum_{B \in \c B_{h_k}(A;k)} \( GT_{\max\set{2h_k+1, \jbar+1}}(B;k) + \sum_{\jbar < j \leq 2h_k} \ol G_j(B; h_k) 2^{-j}\). \label{eq:trailing-upper-bound}
\end{align}
Our lower estimate to $\sigma(A;k)$ is
\begin{equation}
  \label{eq:sigma-lower}
\sigma_-(A;k) = \sigma_+(A;k) - \sum_{B \in \c B_{h_k}(A;k)} \( GT_{\max\set{2h_k+1, \jbar+1}}(A;k) + \sum_{\jbar < j \leq 2h_k} \ol G_j(B; h_k)   2^{-j}\).
\end{equation}
Then $\sigma_-(A;k) \leq \sigma(A;k)$. Note that the computation of $\sigma_-(A;k)$ does not involve computation any terms $G_j(A;k)$ other than the ones used while computing $\sigma_+(A;k)$. However, we do perform a complete search to determine each $\c B_{h_k}(A;k)$, though this is much faster than computing additional $G_j(A;k)$ terms exactly in order to obtain bounds of the same quality.

The strength of Lemma~\ref{lem:error3} lies in that observation that Lemma~\ref{lem:error2} only takes into account the restriction that the last element is \emph{not} in the sum set, whereas Lemma~\ref{lem:error3} additionally takes into account the restriction that the first few elements after $k$ \emph{are} in the sum set.

\subsection{Estimating $\rho$}

Now that we know how to estimate $\sigma(A; k)$ for any particular $(A; k)$, we can obtain the estimates for $\rho(A, B; k) = \sigma(A; k) \sigma(B; k)$ (Proposition \ref{prop:product}) by
\[
	\rho_-(A, B; k) = \sigma_-(A; k) \sigma_-(B; k), \qquad
	\rho_+(A, B; k) = \sigma_+(A; k) \sigma_+(B; k),
\]
where the formulas for $\sigma_+$ and $\sigma_-$ are found in \eqref{eq:sigma-upper} and \eqref{eq:sigma-lower} respectively.
Then, using \eqref{eq:almost-all-bounds}, Proposition \ref{prop:almost-all}, and Lemma \ref{lem:p*=p}, we can obtain the following estimates for $\rho$:
\begin{equation} \label{eq:rho-estimate}
	\sum_{\substack{(A, B; k) \\ k \leq \kbar}} \rho_-(A, B; k)
	\leq \rho
	\leq \(\sum_{\substack{(A, B; k) \\ k \leq \kbar}} \rho_+(A, B; k) \) +  \frac{(3/4)^{\kbar/2}}{1 - \frac{\sqrt 3}2} + 2 \(\frac 3 4\)^{\kbar}
\end{equation}
where the sum is taken over all minimal MSTD fringe pairs $(A, B; k)$ with $k \leq \kbar$.

This completes the description of the algorithm used to estimate $\rho$.

\subsection{Numerical results and comments}

The program was written in Java. All source code are available online at \begin{center}\url{http://web.mit.edu/yufeiz/www/mstd_density_code.zip}\end{center}
 All calculations were performed using exact rational arithmetic. We ran the computation with the following parameters: 
\[
\kbar = 20, \quad 
\jbar = 37, \quad
h_k = \begin{cases} 30, & \text{if } k \leq 10, \\
                               k+10, & \text{if } k > 10. \end{cases}
\]
The entire computation took a combine processing time of approximately one week on a single $2.8$~GHz processor. The results of the computation are shown in Table~\ref{tab:results}. 

% \begin{table}[ht] \centering
% \caption{The number of minimal MSTD fringe pairs.\label{tab:fringe}}

% \begin{tabular}{cc}
% \toprule
% $k$ & \parbox{2in}{\centering The number of minimal MSTD fringe pairs of order $k$} \\ \midrule
% 6 & 8 \\ 
% 7 & 10 \\ 
% 8 & 54 \\ 
% 9 & 106 \\ 
% 10 & 396 \\ 
% 11 & 1034 \\
% 12 & 3120 \\ 
% 13 & 8316 \\ 
% 14 & 26390 \\ 
% 15 & 71594 \\ 
% 16 & 211356 \\ 
% 17 & 612824 \\
% 18 & 1746622 \\
% 19 & 5331566 \\
% 20 & 14747652 \\  \bottomrule
% \end{tabular}
% \end{table}

\begin{table}[ht] \centering
\caption{Results of the computation. The column $\abs{\set{(*,*;k)}}$ contains the number of minimal MSTD fringe pairs of order $k$. The column $\sum \rho_-(*, *; k)$ contains the sum of lower bounds $\rho_-(A, B; k)$ over all minimal MSTD fringe pairs $(A, B)$ of a fixed order $k$, and similarly with the column $\sum \rho_+(*, *; k)$. \label{tab:results}}
\[
\begin{array}{crcc}
  \toprule
  k & \abs{\set{(*,*;k)}} & \sum \rho_-(*, *; k) & \sum \rho_+(*, *; k) \\ \midrule
  6	& 8 	&	0.92959 \x 10^{-4} &	0.93665 \x 10^{-4} \\
  7	& 10 &	0.19475 \x 10^{-4} &	0.19630 \x 10^{-4} \\
  8	& 54 &	0.68801 \x 10^{-4} &	0.69411 \x 10^{-4} \\
  9	& 106 &	0.30178 \x 10^{-4} &	0.30468 \x 10^{-4} \\
  10	& 396 &	0.41411 \x 10^{-4} &	0.41840 \x 10^{-4} \\
  11	& 1034 &	0.34795 \x 10^{-4} &	0.35339 \x 10^{-4} \\ 
  12	& 3120 &	0.29209 \x 10^{-4} &	0.29707 \x 10^{-4} \\
  13	& 8316 &	0.24097 \x 10^{-4} &	0.24529 \x 10^{-4} \\
  14	& 26390 &	0.21456 \x 10^{-4} &	0.21867 \x 10^{-4} \\
  15	&71594 &	0.18176 \x 10^{-4} &	0.18538 \x 10^{-4} \\
  16	& 211356 &	0.13581 \x 10^{-4} &	0.13878 \x 10^{-4} \\
  17	& 612824 &	0.12414 \x 10^{-4} &	0.12701 \x 10^{-4} \\
  18	& 1746622 &	0.08570 \x 10^{-4} &	0.08792 \x 10^{-4} \\
  19	& 5331566 &	0.08035 \x 10^{-4} &	0.08280 \x 10^{-4} \\
  20	& 14747652 &	0.05438 \x 10^{-4} &	0.05624 \x 10^{-4} \\ \midrule
  \Sigma	& &	4.28602 \x 10^{-4} &	4.34262 \x 10^{-4} \\ \bottomrule
\end{array}
\]
\end{table}

Using \eqref{eq:rho-estimate} and the data in Table \ref{tab:results} we obtain
\[
	\rho > 4.286 \x 10^{-4}.
\]
Unfortunately the upper bound that we obtain is rather disappointing, since the error term in the upper estimate in \eqref{eq:rho-estimate} decreases very slowly with $\kbar$:
\[
	\rho < 4.343 \x 10^{-4} + \frac{(3/4)^{20/2}}{1 - \frac{\sqrt 3}2} + 2 \(\frac 3 4\)^{20} < 0.43.
\]

From Monte-Carlo experiments, we know that $\rho$ should be around $4.5\x10^{-4}$, so we see that the weakness in our estimates is in the upper error term as opposed to the sum itself. If we increase $\kbar$, then we should be able to get a better lower bound, but the upper bound would still be far off. The rightmost column sum in Table~\ref{tab:results} represents an upper bound to the best possible lower bound to $\rho$ that we could obtain without increasing $\kbar$. Unfortunately, each increment in $\kbar$ would increase the total computation time by a factor of about four (mostly to due to the search for minimal MSTD fringe pairs). Most of our computation time is spent on complete searches through all subsets of a set (in computing the fringe pairs, $G_j(A;k)$, and $\c B_h(A;k)$), so perhaps it is worthwhile to come up with more efficient search algorithms.

%%%%%%%%%%%%%%%%%%%%%%%%%%%%%%%%%%%%%%%%%%%%%%%%%%%%%%%%%%%%%%%%%%%%%%%%%%%%%%%%%%%%%%%%%%%%%%%%

\section{Extensions to other sum-difference characterizations} \label{sec:extensions}

We have just studied the probability that a random subset $S \subseteq [0,n]$ is an MSTD set. What if we ask finer questions, such as what is the probability that $\abs{S+S} - \abs{S-S} = x$, where $x$ is some fixed integer? Or what is the probability that $S$ is missing exactly $s$ sums and $d$ differences? It turns out that our methods can easily be adapted to deal with all these questions.

Recall from the introduction that
\[	
	\lambda(S) = \( 2n+1 - \abs{S+S}, 2n+1 - \abs{S-S} \)
\]
is the pair consisting of the number of missing sums and the number of missing differences. Fix a subset $\Lambda \subseteq \ZZ_{\geq 0} \x \ZZ_{\geq 0}$. We are interested in the collection
\[
	\{ S \subseteq [0, n] : \lambda(S) \in \Lambda \}.
\]
Let $\rho_n^\Lambda$ be the probability that a uniform random subset $S\subseteq [0, n]$ falls into this collection. In this section we prove Theorem~\ref{thm:gen} showing that $\rho_n^\Lambda$ approaches a limit $\rho^\Lambda$ as $n \to \infty$. By choosing $\Lambda = \{(s,d) : s < d\}$ we get the MSTD problem.

Most of the main ideas for the MSTD case carry over to the general case, so we just sketch the modifications. As with the MSTD problem, we also have a deterministic algorithm for computing arbitrarily good bounds for each limit, though we will not discuss in too much detail the computational aspect as it is similar to Section \ref{sec:computation}. However, even in the case $\Lambda = \{(s,d) : s < d\}$, the general algorithm to be described is much slower than the more specialized algorithm for MSTD sets given earlier. Unlike in Section \ref{sec:computation}, we do not actually carry out the computations, so we make no effort in optimization.

The main difference between the solution of the MSTD case presented earlier and the solution to the general case is that we need to consider a more restrictive analogue of rich sets.

\begin{definition}
Let $k$ and $n$ be positive integers with $2k < n$. Let $S$ be a subset of $[0, n]$. We say that $S$ is \emph{$k$-affluent} if $[k+1, 2n-k-1] \subseteq S+S$ and $[-2n+k+1, 2n-k-1] \subseteq S - S$.
\end{definition}

Whereas rich sets have all the middle sums present, affluent sets additionally have all the middle differences present.

\subsection{Affluent sets with given fringe pair}

In this section we consider the probability that a random $S \subseteq [0, n]$ has a particular fringe profile and is also affluent. The ideas here are very similar to the ones in Sections \ref{sec:fringe} and \ref{sec:semi-rich}. The main difference is that we no longer have the analogue of semi-rich sets since the constraint of being affluent cannot be easily divided into two nearly independent halves.

We need a more general notion fringe pairs to work with affluent sets.

\begin{definition}
A \emph{fringe pair} of order $k$ is a pair of subsets $(A, B)$ of $[0, k]$ (also denoted $(A, B; k)$). We impose the following partial order on fringe pairs: $(A, B; k) > (A', B'; k')$ if $k > k'$ and
\[
	A' = A \cap [0,k'], \quad
	B' = A \cap [0, k'], \quad
	[k'+1, k] \subseteq A+A, B+B, A+B.
\]
\end{definition}

Note that unlike MSTD fringe pairs, we do not require $0 \in A$ or $0 \in B$ here. We previously imposed this requirement as a computational optimization.

We say that a $k$-affluent subset $S\subseteq [0, n]$ has fringe pair $(A, B; k)$ (note that it's the same $k$) if $S \cap [0, k] = A$ and $(n-S)\cap [0, k] = B$.

The partial order for fringe pairs is stronger than the version used to study MSTD sets. As with MSTD fringe pairs, we can speak of \emph{minimal fringe pairs}, as well as the minimal fringe pair of an affluent set. We will count affluent sets by minimal fringe pairs in the same way as we counted rich MSTD sets by minimal MSTD fringe pairs.

Let $(A, B; k)$ be a fringe pair and let $n > 2k$. Let
\[
	\mu_n(A, B; k) = 2^{-n-1} \abs{\set{ S \subseteq [0, n] : S \cap [0, k] = A, (n - S) \cap [0, k] = B, \text{ and } S \text{ is } k \text{-affluent} }}.
\]
Then $\mu_n(A,B;k)$ is the probability that a uniformly random $S \subseteq [0,n]$ (no longer imposing that $0,n \in S$) is $k$-affluent with fringe pair $(A,B;k)$. Let
\[
	\mu(A, B; k) = \lim_{n\to \infty} \mu_n(A, B; k).
\]
The following proposition shows that the limit exists. The result is the analogue of Propositions \ref{prop:sigma-exists} and \ref{prop:product}.

\begin{proposition} \label{prop:mu}
For every $A, B \subseteq [0, k]$, the limit $\mu(A, B; k) = \lim_{n\to \infty} \mu_n(A, B; k)$ exists.
\end{proposition}

\begin{proof}
Assume throughout that $n > 2k$ and $S$ is a uniform random subset of $[0, n]$. We say that $S$ is \emph{$k$-quasi-affluent} if
\begin{align*}
	[k+1, 2n-k-1] \setminus \left[ \floor{\frac n 2}, 2n - \floor{\frac n 2} \right] &\subseteq S+S, \\
	\text{and} \quad [-2n+k+1, 2n-k-1] \setminus \left[-n + \floor{\frac n 2}, n - \floor{\frac n 2} \right] &\subseteq S - S.
\end{align*}
Let $\mu'_n(A, B; k)$ denote the probability that $S$ is $k$-quasi-affluent with fringe pair $(A, B; k)$. If $S$ is $k$-quasi-affluent but not $k$-affluent, then it is necessarily missing some middle sum or middle difference, so we can use Lemma \ref{lem:missing-middle} or an argument analogous to the proof of Proposition \ref{prop:product} to see that this probability goes to zero as $n \to \infty$. In other words,
\[
	\lim_{n \to \infty} ( \mu'_n(A, B; k) - \mu_n(A, B; k) )  = 0.
\]
Thus it suffices to evaluate $\lim_{n \to \infty} \mu'_n(A, B; k)$. Let $m = \floor{\frac n 2} - 1$,
\[
	L = S \cap [0, m], \quad \text{and} \quad
	R = (n - S) \cap [0, m].
\]
Then the condition that $S$ is $k$-quasi-affluent with fringe pair $(A, B; k)$ is equivalent to
\begin{equation} \label{eq:quasi-rich-pair}
	L \cap [0, k] = A, \quad R \cap[0, k] = B, \quad [k+1, m] \subseteq L+L, R+R, L+R.
\end{equation}
So the number of pairs of subsets $(L, R)$ of $[0, m]$ satisfying \eqref{eq:quasi-rich-pair} equals to $\mu'_n(A, B; k)2^{2(m+1)}$.

%If $L$ and $R$ are independent uniform random subsets of $[0, m]$ containing $0$, then $\mu'_n$ is the 
%
%Let $M''_m(A,B; k)$ denote the number of pairs of subsets $(L, R)$ of $[0, m]$ satisfying \eqref{eq:quasi-rich-pair}. Then
%and so
%\[
%	 M'_n(A, B; k)2^{-n+1} = 2^{-2m} M''_m(A, B; k).
%\]

As in the arguments in Section \ref{sec:semi-rich}, we can compute $\mu'_n(A, B; k)2^{2(m+1)}$ by considering the complement to the set of pairs $(L, R)$ satisfying \eqref{eq:quasi-rich-pair}. The complement can be partitioned by the smallest element greater than $k$ missing from any of $L+L, R+R, L+R$. For $j > k$, let $N_j(A, B; k)$ denote the number of pairs $(U, V)$ of $[0, j]$ such that
\begin{gather*} \label{eq:quasi-rich-pair-avoid}
	U \cap [0, k] = A, \quad V \cap[0, k] = B, \quad
	[k+1, j-1] \subseteq L+L, R+R, L+R, \\
	\text{and at least one of } L+L, R+R, L+R \text{ is missing } j.
\end{gather*}
Then
\[
	\mu'_n(A, B; k)2^{2(m+1)} = 2^{2(m-k)} - \sum_{j=k+1}^m N_j(A, B; k) 2^{2(m-j)},
\]
hence
\begin{equation} \label{eq:quasi-rich-compl}
	\mu'_n(A, B; k) = 2^{-2k} - \sum_{j=k+1}^{\floor{n/2} - 1} N_j(A, B; k) 2^{-2j}.
\end{equation}
Since the quantities $\mu'_n(A, B; k)$ and $N_j(A, B; k)$ are all nonnegative, letting $n\to \infty$ shows that the limit
\begin{equation} \label{eq:quasi-rich-limit}
	\mu(A, B; k) =
	\lim_{n \to \infty} \mu'_n(A, B; k) = 2^{-2k} - \sum_{j= k+1}^\infty N_j(A, B; k) 2^{-2j}
\end{equation}
exists.
%The proof that $\mu(A, B; k) > 0$ is very similar to the arguments used in \cite{MO}, so we only sketch the idea. Basically, if we choose a sufficiently large $\ell$ (depending on $k$) and require that $[k+1, \ell] \cup [n-\ell, n-k-1] \subseteq S$, and choose $S \cap [\ell+1, n-\ell-1]$ randomly, then there is a positive probability that all the middle sums and differences are present, thereby making $S$ a $k$-affluent set.
\end{proof}

%Unlike the case of Proposition \ref{prop:sigma-exists}, it is not true that $\mu(A, B; k) > 0$ as we no longer require that $A$ and $B$ contain the element $0$. Indeed, $\mu(\emptyset, \emptyset, 0) = 0$.

Each $\mu(A, B; k)$ can be computed up arbitrary precision using \eqref{eq:quasi-rich-limit}. Indeed, any individual term $N_j(A, B; k)$ can be computed explicitly using a complete search. The tail sum can be bounded using methods analogous to the ones in Section \ref{sec:estimate-sigma}.

\subsection{Almost all sets are affluent}

%As with the MSTD case, we start by restricting to subsets $S \subseteq [0,n]$ containing $0$ and $n$. Let
%\[
%	\rho^\Lambda_{*n} = 2^{-n+1} \abs{\set{ S \subseteq [0,n] : 0, n \in S, \lambda(S) \in \Lambda}}}.
%\]
%The proof of Lemma \ref{lem:p*=p} applies \emph{mutatis mutandis} to the general case. That is, once we show that the limit
%\[
%	\rho^\Lambda_* = \lim_{n \to \infty} \rho^\Lambda_{*n}
%\]
%exists, we then know that the $\rho^\Lambda =  \lim_{n \to \infty} \rho^\Lambda_{n}$ exists and equals $\rho^\Lambda_*$.
Let $\Lambda \subseteq \ZZ_{\geq 0} \x \ZZ_{\geq 0}$ and
\[
	\rho^\Lambda_n = 2^{-n-1} \abs{\set{ S \subseteq [0,n] : \lambda(S) \in \Lambda }}.
\]
For a fringe pair $(A, B; k)$, define
\[
	\lambda(A, B; k) = \(2(k+1) - \abs{(A+A)\cap[0, k]} - \abs{(B+B)\cap[0, k]}, 2(k+1-\abs{(A+B)\cap[0, k]}) \).
\]
It is easy to see that if $S$ is $k$-affluent with fringe pair $(A, B; k)$ then $\lambda(S) = \lambda(A, B; k)$. The following result is the generalization of Propositions \ref{prop:almost-all} and Proposition \ref{prop:formula}.

%Let $\c F(\Lambda)$ denote the set of all minimal fringe pairs $(A, B; k)$ with
%\[
%	(2(k+1) - \abs{(A+A)\cap[0, k]} - \abs{(B+B)\cap[0, k]}, 2(k+1-\abs{(A+B)\cap[0, k]}) \in \Lambda.
%\]
%It is easy to see that if $S$ is an affluent sets with fringe pair in $\Lambda$, then $(2n+1- \abs{S+S},2n+1-\abs{S-S}) \in \Lambda$. 

\begin{proposition} \label{prop:rho-lambda}
As $n\to \infty$, $\rho_n^\Lambda$ converges to a limit $\rho^\Lambda$, and
\begin{equation}\label{eq:rho-lambda}
	\rho^\Lambda = \sum_{\lambda(A, B; k) \in \Lambda}  \mu(A, B; k)
\end{equation}
where the sum is taken over all minimal fringe pairs $(A, B; k)$ satisfying $\lambda(A, B; k) \in \Lambda$.
\end{proposition}

\begin{proof}
%Let $\c F_{\leq \kbar}(\Lambda)$ denote the subset of $\c F(\Lambda)$ consisting of all elements $(A, B; k)$ with $k \leq \kbar$. 
An argument similar to the proof of Proposition \ref{prop:almost-all} shows that
\begin{multline} \label{eq:affluent-bound}
	\sum_{\substack{\lambda(A, B; k) \in \Lambda \\ k \leq \kbar}}  \mu_n(A, B; k)
	\leq 
	\rho^\Lambda_{*n} \\
	\leq
	\sum_{\substack{\lambda(A, B; k) \in \Lambda \\ k \leq \kbar}}  \mu_n(A, B; k) 
		+ \frac{3(3/4)^{\kbar/2}}{2 - \sqrt 3} + 8 \( \frac{3}{4}\)^{\kbar+2} + (n+1) \(\frac{3}{4}\)^{(n-1)/3}.
\end{multline}
The error term on the upper bound uses Lemma \ref{lem:not-affluent}.
Letting $n \to \infty$, and then $\kbar \to \infty$ shows that the limit $\rho^\Lambda_* = \rho^\Lambda$ exists and is equal to the expression in \eqref{eq:rho-lambda}.
\end{proof}

If we want to compute lower and upper bounds for $\rho^\Lambda$, we just let $n \to \infty$ in \eqref{eq:affluent-bound} to get
\[
	\sum_{\substack{\lambda(A, B; k) \in \Lambda \\ k \leq \kbar}}  \mu(A, B; k)
	\leq 
	\rho^\Lambda
	\leq
	\sum_{\substack{\lambda(A, B; k) \in \Lambda \\ k \leq \kbar}}  \mu(A, B; k) 
	 + \frac{3(3/4)^{\kbar/2}}{2 - \sqrt 3} + 8 \( \frac{3}{4}\)^{\kbar+2}.
\]

\begin{proof}[Proof of Theorem \ref{thm:gen}]
The theorem follows almost immediately from Proposition \ref{prop:rho-lambda}. The first assertion is a direct consequence of Proposition \ref{prop:rho-lambda}. The second assertion that $\rho^\Lambda > 0$ as long as $\Lambda$ contains some element $(s, d)$ with $d$ even follows from \cite[Thm.~8]{Hegarty}. For the final assertion, since $\mu(A, B; k) \geq 0$, the sum in \eqref{eq:rho-lambda} can be partitioned by $\lambda(A, B; k)$ to obtain that
\[
	\rho^\Lambda = \sum_{\lambda(A, B; k) \in \Lambda} \mu(A, B; k) 
	= \sum_{(s,d) \in \Lambda} \(\sum_{\lambda(A, B; k) = (s,d)} \mu(A, B; k) \)
	= \sum_{(s,d)\in\Lambda} \rho^{s,d}. \qedhere
\]
\end{proof}

Proposition \ref{prop:rho-lambda} can be used to compute estimates for $\rho^\Lambda$ similar to the MSTD case. The only step that we are missing is bounding the $N_j(A, B; k)$ terms. We omit this discussion since it is very similar to bounding $G_j(A; k)$ as we did in Section \ref{sec:estimate-sigma}.

\section{Structure of a random set characterized by $\Lambda$} \label{sec:middle}

Let $\Lambda \subseteq \ZZ_{\geq 0} \x \ZZ_{\geq 0}$ contain at least one element $(s,d)$ with $d$ even. So $\rho^\Lambda > 0$. In this section, we study the structure of a random subset $S \subseteq [0,n]$ conditioned on $\lambda(S) \in \Lambda$. Our main result, stated below, says that the middle segment of $S$ is nearly unrestricted and independent from the fringe choice. Theorem \ref{thm:middle} formalizes the intuition that the fringe of an MSTD set matters a lot while other elements matter very little.

\begin{theorem} \label{thm:middle}
Let $\Lambda \subseteq \ZZ_{\geq 0} \x \ZZ_{\geq 0}$ where $\Lambda$ contains at least one element $(s,d)$ with $d$ even. Suppose we have an integer sequence $\alpha_n$ satisfying $0 < \alpha_n < n/2$ and $\alpha_n \to \infty$ as $n \to \infty$. 
Let $\epsilon > 0$, then for all sufficiently large $n$ the following is true:

Let $S$ be a uniform random subset of $[0, n]$, $E$ an event that depends only on $S \cap [\alpha_n+1, n - \alpha_n-1]$, and $F$ an event that depends only on $S \cap ([0, \alpha_n] \cup [n - \alpha_n - 1])$. Then
\[
	\abs{ \PP(E \cap F \mid \lambda(S) \in \Lambda ) - \PP(E) \PP(F \mid \lambda(S) \in \Lambda ) }
	\leq (1 + \epsilon) \frac{24(3/4)^{\alpha_n/2}}{(2 - \sqrt 3)\rho^\Lambda}.
\]
\end{theorem}

Note that the bound approaches zero as $n \to \infty$. Intuitively, this says that the structure of the middle portion of a random MSTD set is close to that of an unrestricted set.

\begin{corollary} \label{cor:middle-limit}
Let $\Lambda$ and $\alpha$ satisfy the hypotheses of Theorem \ref{thm:middle}. For each $n$, let $S_n$ be a uniform random subset of $[0, n]$ and $E_n$ an event that depends only on $S_n \cap [\alpha_n+1, n - \alpha_n - 1]$. Suppose that $\lim_{n \to \infty} \PP(E_n)$ exists. Then
\[
	\lim_{n \to \infty} \PP(E_n \mid \lambda_n(S_n) \in \Lambda) = \lim_{n \to \infty} \PP(E_n).
\]
\end{corollary}

\begin{proof}
In Theorem \ref{thm:middle} let $F$ be the event that includes all outcomes.
\end{proof}

In this section we prove Theorem \ref{thm:middle} and give some applications. The proofs are mostly independent of the results in previous sections. Even though we assume the existence of the limit $\rho^\Lambda$, it suffices to know that $\rho_n^\Lambda$ has a positive lower limit. We also use the notion of affluent sets, defined in the beginning of Section~\ref{sec:extensions}.

\subsection{Proof of Theorem \ref{thm:middle}}

We would like to slightly perturb the event on which we are conditioning. The following lemma shows that this modification does not change the probability very much.

\begin{lemma} \label{lem:prob-dep}
Let $A, B, E$ be three events such that $A \subseteq B$ and $\PP(A) > 0$. Then 
\[
	\abs{ \PP(E \mid A) - \PP(E \mid B) } \leq \frac{2 \ \PP(B \setminus A)}{\PP(B)}.
\]
%and $
%Suppose that $\{A_n\}$ and $\{B_n\}$ are two sequences of events such that $A_n \subseteq B_n$ for each $n$, and $\lim_{n\to \infty} \PP(A) = \lim_{n \to\infty} \PP(B) > 0$. If $\{E_n\}$ is another sequence of events, then
%\[
%	\lim_{n \to \infty} \( \PP(E_n \mid A_n) - \PP(E_n \mid B_n) \) = 0.
%\]
\end{lemma}

\begin{proof}
We have
{\allowdisplaybreaks 
\begin{align*} 
	&\abs{ \PP(E \mid A) - \PP(E \mid B) }
	\\&\quad = \abs{ \frac{\PP(E \cap A)}{\PP(A)}  - \frac{\PP(E \cap B)}{\PP(B)} }
	\\&\quad = \frac{\abs{\PP(E \cap A) \PP(B) - \PP(E \cap B)\PP(A)  }}{\PP(A)\PP(B)}
	\\&\quad = \frac{\abs{\PP(E \cap A) \PP(B) - \PP(E \cap A)\PP(A) + \PP(E \cap A)\PP(A) - \PP(E \cap B)\PP(A)  }}{\PP(A)\PP(B)}	
%	\\&\quad \leq \frac{\abs{\PP(E \cap A) - \PP(E \cap B)}\PP(B) + \PP(E \cap B)\abs{\PP(B) - \PP(A)}  }}{\PP(A)\PP(B)}
%	\\&\quad = \frac{\PP(E \cap (B \setminus A))\PP(B) + \PP(E \cap B)\PP(B \setminus A)}{\PP(A)\PP(B)}
%	\\&\quad \leq \frac{\PP(B \setminus A)\PP(B) + \PP(B)\PP(B \setminus A)}{\PP(A)\PP(B)}
	\\&\quad \leq \frac{\PP(E \cap A) \abs{\PP(B) - \PP(A)}+ \abs{\PP(E \cap A) - \PP(E \cap B)}\PP(A) }{\PP(A)\PP(B)}	
	\\&\quad = \frac{\PP(E \cap A) \PP(B \setminus A) + \PP(E \cap (B \setminus A))\PP(A) }{\PP(A)\PP(B)}
	\\&\quad \leq \frac{\PP(A) \PP(B \setminus A) + \PP(B \setminus A)\PP(A) }{\PP(A)\PP(B)}
	\\&\quad \leq \frac{2 \ \PP(B \setminus A)}{\PP(B)},
\end{align*}
}
as desired.
\end{proof}

We would like to slightly perturb the event being conditioned so that it becomes independent of the middle segment of $S$. We do so by adding and removing some non-affluent sets into the event. This is the idea behind the following proposition which leads directly to the theorem.

\begin{proposition} \label{prop:middle1}
Let $\Lambda \subseteq \ZZ_{\geq 0} \x \ZZ_{\geq 0}$, $2k<n$ be positive integers, and $S$ a uniform random subset of $[0, n]$. Assume that $\PP( \lambda(S) \in \Lambda \text{ and } S \text{ is } k\text{-affluent} ) > 0$. Let $E$ be an event that depends only on $S \cap [k+1, n - k-1]$, and $F$ an event that depends only on $S \cap([0, k] \cup [n-k, n])$. Then
\[
	\abs{ \PP(E \cap F \mid \lambda(S) \in \Lambda ) - \PP(E)\PP(F \mid \lambda(S) \in \Lambda) } \leq 
	\frac{8 \ \PP(S \text{ is not } k\text{-affluent})}
		 {\PP( \lambda(S) \in \Lambda \text{ and } S \text{ is } k\text{-affluent} )}.
\]
\end{proposition}

\begin{proof}
Consider the following events:
\begin{align*}
	A &= \{\lambda(S) \in \Lambda \},
\\	B &= \{\lambda(S) \in \Lambda \text{ and } S \text{ is } k\text{-affluent} \},
 \\	C &= \{\exists T \subseteq [0, n], \lambda(T) \in \Lambda, T \text{ is } k\text{-affluent}, \\
 		&\qquad		S \cap [0, k] = T \cap [0, k], S \cap [n-k, n] = T \cap [n -k, n] \},
\\    D &= \{ S \text{ is not $k$-affluent} \}.
\end{align*}
It is easy to see that $B \subseteq A$ and $B \subseteq C$. Furthermore, $A \setminus B \subseteq D$ and $C \setminus B \subseteq D$, the latter follows from the observation that if $C$ occurs and $S$ is $k$-affluent then $S+S = T+T$ and $S-S = T-T$, so that $\lambda(S) = \lambda(T) \in \Lambda$ and hence $B$ occurs as well. 

Applying Lemma \ref{lem:prob-dep} we have
\begin{align*}
	\abs{ \PP(E\cap F \mid A) - \PP(E\cap F \mid B) } &\leq \frac{2 \ \PP(A \setminus B)}{\PP(A)} 
											\leq \frac{2 \ \PP(D)}{\PP(B)}, \\
	\abs{ \PP(E \cap F \mid B) - \PP(E \cap F \mid C) } &\leq \frac{2 \ \PP(C \setminus B)}{\PP(C)}
											\leq \frac{2 \ \PP(D)}{\PP(B)}.
\end{align*}
So combining the two inequalities gives us
\begin{equation} \label{eq:lem-mid1}
	\abs{ \PP(E\cap F \mid A) - \PP(E\cap F \mid C) } \leq \frac{4 \ \PP(D)}{\PP(B)}.
\end{equation}
Similarly, we have
\begin{equation} \label{eq:lem-mid2}
	\abs{ \PP(E)\PP(F \mid A) - \PP(E)\PP(F \mid C) } \leq \frac{4 \ \PP(E) \PP(D)}{\PP(B)}
			\leq \frac{4 \ \PP(D)}{\PP(B)}.
\end{equation}
Now, $E$ depends only on $S \cap [k+1, n - k-1]$, while $F$ and $C$ depend only on $S \cap([0, k] \cup [n-k, n])$. So $E$ is independent from $F \cap C$. Thus $\PP(E\cap F \mid C) = \PP(E)\PP(F \mid C)$. Then combining \eqref{eq:lem-mid1} and \eqref{eq:lem-mid2} gives us
\begin{align*}
\abs{ \PP( E\cap F \mid A) - \PP(E)(F \mid A) }
& \leq \abs{ \PP( E\cap F \mid A) - \PP(E \cap F \mid C) } + \abs{ \PP(E)(F \mid A) - \PP(E) \PP(F \mid C) }
\\ & \leq \frac{8 \ \PP(D)}{\PP(B)},
\end{align*}
as desired.
\end{proof}

\begin{proof}[Proof of Theorem \ref{thm:middle}]
Let $S_n$ denote a uniform random subset of $[0, n]$. Using Proposition \ref{prop:middle1}, it suffices to show that 
\[
	\limsup_{n \to \infty} \frac{8 \ \PP(S_n \text{ is not } \alpha_n\text{-affluent}) (3/4)^{-\alpha_n/2}}
		 {\PP( \lambda(S_n) \in \Lambda \text{ and } S_n \text{ is } \alpha_n\text{-affluent} )} 
		 \leq \frac{24}{(2 - \sqrt 3)\rho^\Lambda}.
\]
By Lemma \ref{lem:not-affluent} we have
\begin{equation} \label{eq:prob-not-affluent}
	\PP(S_n \text{ is not } \alpha_n\text{-affluent})
		\leq \frac{3(3/4)^{\alpha_n/2}}{2 - \sqrt 3} + 8 \( \frac{3}{4}\)^{\alpha_n+2} + (n+1) \(\frac{3}{4}\)^{(n-1)/3},
\end{equation}
so that
\[
	\limsup_{n \to \infty} \PP(S_n \text{ is not } \alpha_n\text{-affluent}) (3/4)^{-\alpha_n/2} \leq  \frac{3}{2 - \sqrt 3}.
\]
By \eqref{eq:prob-not-affluent} and Theorem~\ref{thm:gen} we have
\[
	\lim_{n \to \infty} \PP( \lambda(S_n) \in \Lambda \text{ and } S_n \text{ is } \alpha_n\text{-affluent} )
	= \lim_{n \to \infty} \PP( \lambda(S_n) \in \Lambda )
	= \rho^\Lambda.
\]
The theorem then follows.
\end{proof}

\subsection{Applications}

In this section we explore some applications of Theorem \ref{thm:middle}. 

Miller, Orosz, and Scheinerman \cite{MOS} conjectured that, for a fixed constant $0 < c < 1/2$, and $k_n$ varying with $n$ satisfying $c n < k_n < n - c n$, we have
\[
	\lim_{n \to \infty} \frac{ \abs{\set{S \subseteq [0, n] : k_n \in S \text{ and } S \text { is MSTD} }}}
							 {\abs{\set{S \subseteq [0, n] : S \text { is MSTD} }}} 
						= \frac 1 2.
\]
It was also asked if we could replace the condition $c n \leq k_n \leq n - c n$ by $\alpha_n < k_n < n - \alpha_n$ for some function $\alpha$. The following result answers these questions. Recall that taking $\Lambda = \{(s,d) : s < d\}$ gives us MSTD sets.

\begin{corollary} \label{cor:prob-middle}
Let $\Lambda$ and $\alpha$ satisfy the hypotheses of Theorem \ref{thm:middle}. For each $n$, let $S_n$ be a uniform random subset of $[0, n]$. If $k_n$ is a sequence satisfying $\alpha_n < k_n < n - \alpha_n$, then
\[
	\lim_{n \to \infty} \PP( k_n \in S_n \mid \lambda(S_n) \in \Lambda ) = \frac 1 2.
\]
\end{corollary}

\begin{proof}
In Corollary \ref{cor:middle-limit}, let $E_n$ be the event $\{k_n \in S_n\}$.
\end{proof}

Now we give some results about the size of a random subset $S \subseteq [0,n]$ satisfying $\lambda(S) \in \Lambda$. Because fringe elements do not contribute significantly to $\abs{S}$, our intuition tells us that the size of the set should behave similar to an unrestricted binomial distribution. The next two results confirm this intuition. In the variance part of the next Proposition, we actually need to set the fringe event $F$ in Theorem~\ref{thm:middle} to be something nontrivial, thereby using the full power of the theorem.

\begin{proposition} \label{prop:exp-var}
Let $\Lambda \subseteq \ZZ_{\geq 0} \x \ZZ_{\geq 0}$ contain at least one $(s,d)$ with $d$ even. For each $n$, let $S_n$ be a uniform random subset of $[0,n]$. Then
\begin{equation} \label{eq:exp}
	\EE [ \abs{S_n} \mid \lambda(S_n) \in \Lambda]  = \frac{n+1}{2} + O(\log n)
\end{equation}
and
\begin{equation} \label{eq:var}
	\Var( \abs{S_n} \mid \lambda(S_n) \in \Lambda) = \frac{n+1}{4} + O((\log n)^2)
\end{equation}
where the constants in the big-$O$ may depend on $\Lambda$.
\end{proposition}

\begin{proof}
Choose $\alpha_n = \floor{c\log n}$ for some constant $c > \frac{4}{\log(4/3)}$. Let $S^\alpha_n = S_n \cap [\alpha_n+1, n-\alpha_n-1]$. Applying Theorem \ref{thm:middle} to the events $E = \{ k_n \in S_n\}$ and $F$ the event of all outcomes, we get
\begin{align*}
	\abs{ \EE [ \abs{S^\alpha_n} \mid \lambda(S_n) \in \Lambda] - \frac{n-1-2\alpha_n}2 }
&	\leq \sum_{k = \alpha_n+1}^{n-\alpha_n-1} \abs{ \PP \( k \in S_n \mid \lambda(S_n) \in \Lambda \) - \frac 1 2}
\\&	= O\(n (3/4)^{\alpha_n/2} \)
\\&	\to 0, \text{ as } n \to \infty.
\end{align*}
Thus
\begin{align*}
	\EE [ \abs{S_n} \mid \lambda(S_n) \in \Lambda] 
&	= \EE [ \abs{S_n^\alpha} \mid \lambda(S_n) \in \Lambda] + \EE [ \abs{S_n \setminus S^\alpha_n} \mid \lambda(S_n) \in \Lambda]
\\&	= \frac{n-1 - 2\alpha_n}2 + o(1) + O(\alpha_n)
\\&	= \frac {n+1} 2 + O(\log n).
\end{align*}
This proves \eqref{eq:exp}.

Next, for the variance, we have
\begin{align*}
	\Var( \abs{S_n} \mid \lambda(S_n) \in \Lambda)
&	= \EE \left[ \(\abs{S_n} - \EE[\abs{S_n}] \)^2 \mid \lambda(S) \in \Lambda \right]
\\&	= \EE \left[ \(\abs{S_n} - \tfrac {n+1} 2 + O(\log n) \)^2 \mid \lambda(S) \in \Lambda \right]
\\&	= \EE \left[ \(\abs{S_n} - \tfrac {n+1} 2 \)^2 \mid \lambda(S) \in \Lambda \right] 
\\&\qquad		+ O(\log n) \EE \left[ \(\abs{S_n} - \tfrac {n+1} 2 \) \mid \lambda(S) \in \Lambda \right]
			+ O((\log n)^2)
\\&	= \EE \left[ \(\abs{S_n} - \tfrac {n+1} 2 \)^2 \mid \lambda(S) \in \Lambda \right] + O((\log n)^2).
\end{align*}
For each $i \in [0, n]$, let $X_i$ be the indicator random variable which is 1 if $i \in S$ and $0$ otherwise. Then
\begin{align}
	\EE \left[  \(\abs{S_n} - \tfrac {n+1} 2 \)^2 \big\vert\ \lambda(S) \in \Lambda \right]
&	= \EE \left[ \(\sum_{i=0}^n \(X_i - \tfrac12 \) \)^2 \Bigg\vert\ \lambda(S) \in \Lambda \right] \nonumber
\\&	= \sum_{i=0}^n \sum_{j=0}^n \EE \left[ \(X_i - \tfrac12 \)\(X_j - \tfrac12 \)  \mid \lambda(S) \in \Lambda \right]. \label{eq:var-sum}
\end{align}

Next we analyze each term $\EE \left[ \(X_i - \tfrac12 \)\(X_j - \tfrac12 \)  \mid \lambda(S) \in \Lambda \right]$ using Theorem \ref{thm:middle}. There are several cases to consider.

Suppose that $i, j \in [\alpha_n+1, n - \alpha_n-1]$. For any event $E$ that depends on $S \cap \{i, j\}$, we have
\[
	\abs{ \PP(E \mid \lambda(S) \in \Lambda ) - \PP(E) } = O((3/4)^{\alpha_n/2}).
\]
Thus,
\begin{align*}
	\EE \left[ \(X_i - \tfrac12 \)\(X_j - \tfrac12 \)  \mid \lambda(S) \in \Lambda \right]
&	= \EE \left[ \(X_i - \tfrac12 \)\(X_j - \tfrac12 \) \right] + O((3/4)^{\alpha_n/2})
\\&	= O((3/4)^{\alpha_n/2}) + \begin{cases}
			\tfrac14 , & \text{if } i = j, \\
			0, & \text{if } i \neq j. \end{cases}
\end{align*}

Next, suppose that $i \in [\alpha_n+1, n - \alpha_n-1]$ and $j \notin [\alpha_n+1, n - \alpha_n-1]$ (or vice-versa). If event $E$ is either $\{ i \in S \}$ or $\{i \notin S\}$ and event $F$ is either $\{j \in S\}$ or $\{j \notin S\}$, then
\[
	\abs{ \PP(E\cap F \mid \lambda(S) \in \Lambda ) - \PP(E)\PP(F \mid \lambda(S) \in \Lambda ) } = O((3/4)^{\alpha_n/2}).
\]
Also
\begin{align*}
	\EE \left[ \(X_i - \tfrac12 \)\(X_j - \tfrac12 \)  \mid \lambda(S) \in \Lambda \right]
&	= \EE \left[ X_i - \tfrac12 \right] \EE\left[ X_j - \tfrac12  \mid \lambda(S) \in \Lambda \right] + O((3/4)^{\alpha_n/2})
\\&	= O((3/4)^{\alpha_n/2}).
\end{align*}

Finally, if $i, j \notin [\alpha_n+1, n - \alpha_n-1]$ then we simply use the crude approximation
\[
	- \tfrac14 \leq \EE \left[ \(X_i - \tfrac12 \)\(X_j - \tfrac12 \)  \mid \lambda(S) \in \Lambda \right] \leq \tfrac14.
\]

Combining all three cases and continuing \eqref{eq:var-sum} we get
\begin{align*}
	\Var( \abs{S_n} \mid \lambda(S_n) \in \Lambda)
&	= \EE \left[ \(\abs{S_n} - \tfrac {n+1} 2 \)^2 \big\vert\ \lambda(S) \in \Lambda \right] + O((\log n)^2)
\\& = \sum_{i=0}^n \sum_{j=0}^n \EE \left[ \(X_i - \tfrac12 \)\(X_j - \tfrac12 \)  \mid \lambda(S) \in \Lambda \right] + O((\log n)^2)
\\&	= \frac{n+1}4 + O(n^2 (3/4)^{\alpha_n/2}) + O(\alpha_n^2) + O((\log n)^2)
\\&	= \frac{n+1}4 + O((\log n)^2). \tag*{\qedhere}
\end{align*}
\end{proof}

The next result shows that the size of $S$ follows a central limit theorem.

\begin{proposition} \label{eq:CLT}
Let $\Lambda \subseteq \ZZ_{\geq 0} \x \ZZ_{\geq 0}$ contain at least one $(s,d)$ with $d$ even. For each $n$, let $S_n$ be a uniform random subset of $[0,n]$. Then, for any real number $t$, we have
\[
	\lim_{n \to \infty} \PP \( \abs{S_n} < \frac{n + t\sqrt{n}}2 \ \bigg\vert\ \lambda(S_n) \in \Lambda  \)
	 = \Phi(t)
\]
where $\Phi(t)$ is the standard normal distribution.
\end{proposition}

\begin{proof}
Choose any $\alpha_n = o(\sqrt n)$ with $\alpha_n \to \infty$. Let $S^\alpha_n$ denote $S_n \cap [\alpha_n+1, n-\alpha_n-1]$. We have
\begin{align*}
	\PP \( \abs{S^\alpha_n} < \frac{n + t\sqrt{n}}2 - 2\alpha_n - 2 \ \bigg\vert\ \lambda(S_n) \in \Lambda  \)
&	\leq 	\PP \( \abs{S_n} < \frac{n + t\sqrt{n}}2 \ \bigg\vert\ \lambda(S_n) \in \Lambda  \)
\\&	\leq	\PP \( \abs{S^\alpha_n} < \frac{n + t\sqrt{n}}2 \ \bigg\vert\ \lambda(S_n) \in \Lambda  \).
\end{align*}
Using Corollary \ref{cor:middle-limit} and the Central Limit Theorem, we find that
\[
	\lim_{n \to \infty} \PP \( \abs{S^\alpha_n} < \frac{n + t\sqrt{n}}2 \ \bigg\vert\ \lambda(S_n) \in \Lambda  \)
	= \lim_{n \to \infty} \PP \( \abs{S^\alpha_n} < \frac{n + t\sqrt{n}}2 \)
	= \Phi(t).
\]
Similarly,
\[
	\lim_{n \to \infty} \PP \( \abs{S^\alpha_n} < \frac{n + t\sqrt{n}}2 - 2\alpha_n - 2 \ \bigg\vert\ \lambda(S_n) \in \Lambda  \)
	= \Phi(t).
\]
Therefore
\[
	\PP \( \abs{S_n} < \frac{n + t\sqrt{n}}2 \ \bigg\vert\ \lambda(S_n) \in \Lambda  \) = \Phi(t). \qedhere
\]
\end{proof}

%%%%%%%%%%%%%%%%%%%%%%%%%%%%%%%%%%%%%%%%%%%%%%%%%%%%%%%%%%%%%%%%%%%%%%%%%%%%%%%%%%%%%%%%%%%%%%%%%%

\section{Conclusion and discussion} \label{sec:conclusion}

This paper explores the intuition about the structure of a random MSTD set, namely that its fringe elements are significant while its middle elements are not. Consequently, we can compute the proportion of MSTD sets by searching through all desirable fringe pairs and then sum up the contributions from each fringe pair. We were also able to make some precise statements about how the middle elements are nearly unrestricted and independent from the fringe elements.

More generally, our results apply to any characterization $\Lambda$ on the number of missing sums and the number of missing differences of $S \subseteq \{0, 1, \dots, n\}$. Our methods can also be modified to deal with the following two extensions, though we choose not to discuss them in order to keep the arguments simple.
\begin{itemize}
	\item Our paper is based on the model where each element of $\{0, 1, \dots, n\}$ is chosen independently with probability $1/2$. Our results can be modified to deal with the model where the probability is some other constant (independent of $n$).
	\item We can place additional constraints on the fringe of $S$. For example, in addition to requiring $\lambda(S) \in \Lambda$, we may further require that $0, 1, n \in S$ and $4, n-1 \notin S$. This amounts to including or excluding a certain subset of prefix-suffix pairs.
\end{itemize}

Our method currently does not easily extend to the model where the each element is chosen with probability $p(n)$ varying with $n$. For results in this direction, Hegarty and Miller \cite{HM} showed that if $p(n) \to 0$ and $n^{-1} = o(p(n))$, then a random subset almost always has more sums than differences. It would be interesting to see if there are any analogues of Theorem \ref{thm:middle} other than in the uniform model with constant probability.

We showed that each limit $\rho^\Lambda$ can be computed deterministically up to arbitrary precision. However, in practice, the convergence is very slow since each term requires a complete search. Also error bounds such as Lemma \ref{lem:missing-middle} are too weak to give good numerical results. In the MSTD case we were able to substantially speed up the computation by splitting a rich set into two semi-rich sets and then analyzing each half separately. Unfortunately, in the general case, there does not seem to be a good way to split up an affluent set. Consequently, we expect the computation in the general case to be much slower.

It would nice to find some optimization that could substantially speed up the algorithm. For instance, perhaps we do not have to perform so many complete searches, or perhaps there is some way to divide an affluent set into nearly independent parts. It would also be nice to have a tigher upper bound than what is provided by Lemma \ref{lem:missing-middle}.

In practice, if we wish to estimate any $\rho^\Lambda$, the easiest and quickest way would be to run a Monte Carlo simulation. However, this has the disadvantage of not being able to give any provable bounds.

We conclude with some possible further questions.
\begin{enumerate}
	\item For each fringe pair $(A, B; k)$, can we give an explicit construction of a family of rich/affluent sets that occupy $\Omega(1)$ density?
	\item What can we say if we choose to characterize $S$ by $(\abs{S+S}, \abs{S-S})$ instead of the number of missing sums and differences? In this case, which subsets of $\ZZ_{\geq 0} \x \ZZ_{\geq 0}$ give interesting results?
	\item How quickly does $\rho^\Lambda_n$ converge to $\rho^\Lambda$? Our proofs do not say anything about this. The convergence mentioned in this paper is the convergence of the computed numerical bound, which depends on the order $k$  of the fringe pairs as opposed to $n$.
	\item For which $\Lambda$ is the sequence $\{\rho_n^\Lambda\}$ monotonic? Martin and O'Bryant \cite{MO} suggest perhaps it is monotonically increasing for $\{(s,d) : s < d\}$ and $\{(s,d) : s > d\}$, while monotonically decreasing for $\{(s,d) : s= d\}$. Is the sequence $\{\rho_n^\Lambda\}$ always eventually monotonic? When does it approach the limit $\rho^\Lambda$ from above and when does it approach the limit from below?
	\item Can we improve the error term in Proposition \ref{prop:exp-var} for the expectation and variance of $\abs{S}$? For which $\Lambda$ is the error term asymptotically tight?
\end{enumerate}

%%%%%%%%%%%%%%%%%%%%%%%%%%%%%%%%%%%%%%%%%%%%%%%%%%%%%%%%%%%%%%%%%%%%%%%%%%%%%%%%%%%%%%%%%%%%%%%%%%

\section*{Acknowledgments}

This research was carried out at the University of Minnesota Duluth under the supervision of Joseph Gallian with the financial support of the National Science Foundation and the Department of Defense (grant number DMS 0754106), the National Security Agency (grant number H98230-06-1-0013), and the MIT Department of Mathematics. The author would like to thank Joseph Gallian for his encouragement and support. The author would also like to thank Reid Barton and Nathan Kaplan for reading the paper and making valuable suggestions.

%%%%%%%%%%%%%%%%%%%%%%%%%%%%%%%%%%%%%%%%%%%%%%%%%%%%%%%%%%%%%%%%%%%%%%%%%%%%%%%%%%%%%%%%%%%%%%%%

%\bibliographystyle{alpha}
\bibliographystyle{amsplain}
\bibliography{../references}

\end{document}